\documentclass[matann2]{svjour}

\usepackage{amsmath}
\usepackage{amsfonts}
\usepackage{latexsym}
\usepackage{mathrsfs}
\usepackage{times}
\usepackage{array}
\usepackage{amscd}
\usepackage{a4}
\usepackage{amssymb}
\usepackage{pstricks}
\usepackage{amsrefs}

\input{boxedeps} 
\SetEPSFDirectory{} 
\HideDisplacementBoxes
\SetRokickiEPSFSpecial  
\DeclareMathAlphabet{\ams}{U}{msb}{m}{n}
\DeclareMathAlphabet{\goth}{U}{euf}{m}{n}

\def\sl{\text{SL}}

\def\ml{M\kern-1pt L}

\def\fix{\text{Fix}}

\def\id{\text{id}}

\def\im{\text{im}\,}
\def\aut{\text{Aut}}

\def\ov{\overline}

\def\wh{\widehat}

\def\dom{\text{dom}}

\def\II{\mathcal I}
\def\JJ{\mathcal J}
\def\BB{\mathcal B}

\def\FF{\mathcal F}
\def\AA{\mathcal A}

\def\HH{\mathcal H}

\def\PP{\mathcal P}

\def\TT{\mathcal T}

\def\cS{\mathcal S}

\def\SS{\goth{S}}
\def\BBB{\goth{B}}
\def\XXX{\goth{X}}

\def\ve{\varepsilon}

\def\aa{\alpha}

\def\bb{\beta}
\def\ss{\sigma}

\def\ve{\varepsilon}
\def\Om{\Omega}

\def\wh{\widehat}

\def\Z{\ams{Z}}
\def\R{\ams{R}}

\def\T{\ams{T}}

\def\G{\ams{G}}
\def\M{\ams{M}}

\def\vv{\mathbf{v}}
\def\uu{\mathbf{u}}
\def\xx{\mathbf{x}}

\def\quo{/\kern -.45em\sim}

\def\ds{\displaystyle}

\def\Langle{\langle\kern -2pt\langle}
\def\Rangle{\rangle\kern -1.9pt\rangle}

\def\lf{\lfloor}
\def\rf{\rfloor}
\newpsobject{showgrid}{psgrid}{subgriddiv=1,griddots=10,gridlabels=6pt,gridcolor=red}
\title{Partial symmetry, reflection monoids and Coxeter groups}

\author{Brent Everitt and John Fountain
\thanks{The authors are very grateful to the referee, whose careful reading and numerous
suggestions have significantly improved the exposition of this paper.
The first author would also like to thank 
Bob Howlett 
for helpful suggestions.
Some of the results of this paper were 
obtained while he was visiting
the Institute for Geometry and its Applications, 
University of Adelaide, Australia. 
He is grateful for their hospitality. The second author would like to
thank Mark Kambites and Ben Steinberg for some helpful discussions, and Chris Hollings for
his help with several calculations. A grant from the Royal Society made
it possible for him to visit the University of Adelaide to continue the
work reported here. He would like to express his gratitude to the members of
the Glenelg Mathematics Institute for their kindness and hospitality
during his visit to Adelaide.}
}

\institute{
Department of Mathematics, University of York, York
YO10 5DD, United Kingdom. \email{bje1@york.ac.uk} (Brent Everitt), 
\email{jbf1@york.ac.uk} (John Fountain).
}

\titlerunning{}
\authorrunning{Brent Everitt and John Fountain}

\begin{document}

\maketitle


\begin{abstract}
This is the first of a series of papers in which 
we initiate and develop the theory of reflection monoids, 
motivated by the theory of reflection
groups. The main results 
identify a number of important inverse semigroups as reflection monoids,
introduce new examples, and determine their orders.
\end{abstract}


\section*{Introduction}\label{section:introduction}

The symmetric group $\SS_n$ comes in many guises: as the permutation group of the
set $\{1,\ldots,n\}$; as the group
generated by reflections in the hyperplanes $x_i-x_j=0$ of 
an $n$-dimensional Euclidean space;
as the Weyl group 
of the reductive algebraic group $GL_n$, or (semi)simple group $SL_{n+1}$,
or simple Lie algebra $\goth{sl}_{n+1}$; as the Coxeter group
associated to Artin's braid group, $\ldots$

If one thinks of $\SS_X$ as the group of (global) symmetries of $X$, then the partial
symmetries naturally lead one to consider the \emph{symmetric inverse monoid} 
$\II_X$, whose
elements are the partial bijections $Y\rightarrow Z$ $(Y,Z\subset X)$.
It too has many other faces. It arises in its incarnation 
as the ``rook monoid'' as 
the Renner monoid of the reductive algebraic monoid $M_n$.
An 
associated Iwahori theory and representations 
have been worked out by Solomon \cite{Solomon02,Solomon90}. 
There is a braid connection too, with $\II_n$ naturally associated to 
the inverse monoid of ``partial braids'' defined recently in \cite{Easdown04}.

But what is missing is a realization of $\II_n$ as a
``partial reflection monoid'', or indeed, a definition and theory
of partial mirror symmetry and the monoids generated by partial reflections
that generalizes the theory of reflection groups.

Such is the purpose of the present paper. Reflection monoids are defined
as monoids generated by certain partial linear isomorphisms 
$\aa:X\rightarrow Y$ ($X,Y$ subspaces of $V$), 
that are the
restrictions (to $X$) of reflections. Initially one is faced with many 
possibilities, with the challenge being to impose enough structure for
a workable theory while still encompassing as many interesting examples as 
possible. 
It turns out that a solution is to consider monoids of partial linear
isomorphisms whose domains form a $W$-invariant semilattice for
some reflection group $W$ acting on $V$. 

Two pieces of data will thus go into the definition of a reflection monoid:
a reflection \emph{group\/} and a collection of well behaved domain subspaces
(see \S \ref{section:reflectionmonoids} for the precise definitions).
What results is a theory of reflection monoids for which our main
theorems in this paper determine their orders
and identify the natural examples.

For instance, just as $\SS_n$ is the reflection group associated to the 
type $A$ root system, so
now $\II_n$ becomes the reflection monoid associated to the type $A$ root
system, and where the domains
form a $\SS_n$-invariant Boolean lattice (see \S\ref{section:symmetricgroup}). 
This in fact turns out to be a common feature: if the reflection group
is $\SS_n$ and the domain subspaces are allowed to
vary, one often gets families of monoids pre-existing in the
literature. These families are then generalized by replacing $\SS_n$ by an
arbitrary reflection group.
Thus, the group of signed permutations
of $\{1,\ldots,n\}$ is the Weyl group of type $B$, and the inverse monoid
$\II_{\pm n}$ of \emph{partial\/} signed permutations is a reflection 
monoid of type $B$,
with again a Boolean lattice of domain subspaces. 

Another interesting class of examples arises from the theory of hyperplane arrangements.
The 
\emph{reflection\/} or \emph{Coxeter arrangement monoids\/} have as their 
input data a reflection group $W$
and for the domains, the intersection lattice of the reflecting hyperplanes of $W$.
These intersection lattices possess many beautiful combinatorial and 
algebraic properties (see \cite{Orlik92}).
Thus, the Coxeter arrangement monoids tie up reflection groups and the intersection
lattices of their reflecting hyperplanes in one very natural algebraic
object.
Our main result here (Theorem
\ref{section:arrangementmonoids:result250} below) is that the orders
of these 
reflection monoids are the sum of the indices of the parabolic
subgroups of the original reflection group.

By the ``rigidity of tori'', a maximal torus $T$ in a linear algebraic group
$\G$ has automorphisms a finite group, the Weyl group of $\G$, 
and this is a reflection
group in the space $\XXX(T)\otimes\R$, where $\XXX(T)$ is the character group
of the torus. 
For a linear algebraic monoid $\M$, our view is that there are two finite inverse monoids
needed to play an equivalent role. One, the Renner monoid, is already well known
for the important role it plays in the Bruhat decomposition of $\M$. The other is
a reflection monoid
in $\XXX(T)\otimes\R$, where the extra piece of data, the semilattice of domain spaces, 
comes from the character semigroup $\XXX(\ov{T})$ of the Zariski
closure of $T$. Our main result in this direction, Theorem
\ref{section:renner:result100}, is that these two
are pretty  
closely related.


The paper is organized as follows: \S\ref{section:reflectionmonoids} 
contains the basics and our first main
result, on the orders of a large class of finite inverse monoids; 
\S\ref{section:symmetricgroup} has a 
couple of examples based around the symmetric group. The idea is that
``anyone'' could read these sections, irrespective of whether
they have an interest in reflection groups or semigroups. We fix our
notation concerning 
reflection groups in \S\ref{section:reflectiongroups},
and this allows us to 
explore our first two families of examples in \S\S \ref{section:booleanmonoids}
and \ref{section:arrangementmonoids}:
the Boolean and Coxeter arrangement reflection monoids. In 
\S\ref{section:generalities}, we give some fundamental abstract properties. The reflection 
monoid associated to a reductive algebraic monoid and its relation
to the Renner monoid are the subject of \S\ref{section:renner}.
The last section, as the Bourbaki-ism in its title suggests, 
is a portmanteau of results of independent interest.

In the sequel \cite{Everitt:Fountain} to this paper,
a general presentation is derived (among other things) 
for reflection monoids.
This presentation is determined explicitly, and massaged a little more, 
for the Boolean and Coxeter arrangement monoids associated to 
the Weyl groups. The benchmark here is provided by a classical 
presentation \cite{Popova61} for the symmetric inverse monoid $\II_n$, 
which we rederive
in its new guise as the ``Boolean reflection monoid of type $A$''.

\section{Monoids of partial linear isomorphisms and reflection monoids}
\label{section:reflectionmonoids}

The symmetric group $\SS_X$ and the general linear group $GL(V)$ measure the symmetry 
of a set and a vector space. We start this section with two algebraic objects
that measure instead \emph{partial\/} symmetry. One is (reasonably) well
known, the other less so, but nevertheless implicit in the area. 

For a non-empty set $X$, a \emph{partial permutation} is a bijection 
$Y\to Z$, where $Y,Z$ are subsets of $X$.
We allow $Y$ and $Z$ to be empty, so that the empty function 
$0_X:\varnothing\rightarrow\varnothing$ is regarded
as a partial permutation. The set of all partial permutations of $X$ is made
into a monoid with zero $0_X$ using the usual rule for composition of partial
functions: it is called the \emph{symmetric inverse monoid} on $X$, and denoted by $\II_X$.
If $X =\{1,2,\dots,n\}$, we write $\II_n$ for $\II_X$. See \cite{Howie95}*{\S 5}.

Now let $k$ be a field and $V$ a vector space over $k$.
A \emph{partial linear isomorphism\/} of $V$ is a vector space isomorphism 
$Y\rightarrow Z$, where $Y,Z$ are vector subspaces of $V$. 
The set of partial linear isomorphisms of $V$ is also made
into a monoid using composition of partial
functions (and with zero the linear isomorphism $\mathbf{0}\rightarrow\mathbf{0}$, from the zero
subspace to itself). We call it the \emph{general linear monoid\/} on
$V$ and denote it by $\ml(V)$.

It is possible to toggle back and forth between these two monoids, 
using the inclusions $\ml(V)\subset\II_V$ and $\II_X\subset\ml(V)$, for $V$ the 
$k$-space with basis $X$. Either can be taken as the motivating example of
an \emph{inverse monoid\/} \cites{Howie95,Lawson98}: 
a monoid $M$ such that for all $a\in M$
there is a unique $b\in M$ such that $aba = a$ and $bab = b$. The
element $b$ is the \emph{inverse\/} of $a$ and is denoted by $a^{-1}$.
Intuitively, if $a$ is the partial map $Y\rightarrow Z$, then 
$a^{-1}$ is the inverse partial map $Z\rightarrow Y$, and it is precisely
in order to capture this idea of ``local inverses'' that the notion
of inverse monoid was formulated.

When $X$ is finite, or $V$ finite dimensional, any partial permutation/isomorphism
$Y\rightarrow Z$ can be obtained by restricting to $Y$ a full permutation/isomorphism
$g:X\rightarrow X$. We will write $g_Y^{}$ for the partial map
with domain $Y$ and effect that of restricting $g$ to $Y$.
Equivalently, $g_Y^{}=\ve_Y^{}g$ where $\ve:X\rightarrow X$ is the identity
and $\ve_Y^{}:Y\rightarrow Y$ its restriction to $Y$, a
\emph{partial identity\/}. Thus every partial map is the product of an idempotent
and a unit.
One has to be careful with such representations: $g_Y^{}=h_Z^{}$ if and only
if $Y=Z$ and $gh^{-1}$ is in the \emph{isotropy group\/} 
$G_Y=\{g\in G\,|\,\vv g=\vv,\text{ for all }\vv\in Y\}$
of $Y$. We have 
\begin{equation}
  \label{section:reflectionmonoids:equation100}
g_Y^{}h_Z^{} = (gh)_{Y \cap\kern1pt Zg^{-1}}^{},
\end{equation}
and $(g_Y^{})^{-1} =(g^{-1})_{Y\kern-1pt g}^{}$.
From now on, \emph{all\/} our vector spaces will be finite dimensional.

Again, generalities are suggested by these primordial examples: a monoid $M$ is
\emph{factorizable\/} if $M=EG$ with $E$ the idempotents and $G$ the units of $M$.
The role of the isotropy group is played by the idempotent stabilizer 
$G_e=\{g\in G\,|\,eg=e\}$, and we have equality $e_1g_1=e_2g_2$ if and only if
$e_1=e_2$ and $g_2g_1^{-1}\in G_{e_2}$. 
The units act on the idempotents:
if $e\in E$ and $g\in G$ then $g^{-1}eg\in E$ (with
$g^{-1}\ve_Y^{}g=\ve_{Yg}^{}$ in the 
examples above).

Looking a little more closely at the domain $Y \cap\kern1pt Zg^{-1}$ of $g_Y^{}h_Z^{}$
suggests the following:

\begin{definition}\label{section:reflectionmonoids:definition50}
Let $V$ be a vector space and $G\subset GL(V)$ a group. 
A collection $\cS$ of subspaces of $V$ is called a \emph{system in $V$ for $G$\/}
if and only if 
\begin{description}
\item [(S1).] $V\in\cS$,
\item [(S2).] $\cS G=\cS$, ie: $Xg\in\cS$ for
any $X\in\cS$ and $g\in G$, and
\item [(S3).] if $X, Y\in\cS$ then $X\cap Y\in\cS$. 
\end{description}
\end{definition}

If $\cS_i\,(i\in I)$ is a family of systems for $G$ then
$\bigcap\cS_i$ is too, and thus for any set $\Omega$ of subspaces 
we write $\langle\Omega\rangle_G$ for the intersection of all systems
for $G$ containing $\Omega$, and call this the system 
for $G$ {\em generated\/} by $\Omega$.

Clearly one can always find trivial systems for $G$--just take $V$ 
itself for instance--as well as plenty of examples, about which we
can't say a great deal: 
the system $\langle\Omega\rangle_G$ for $G$ generated
by any set $\Omega$ of subspaces.
There is one system though that is intrinsic to $G$,
encoding some of its structure: for $H$ a subgroup of $G$ let 
$\text{Fix}(H)=\{\vv\in V\,|\, \vv g=\vv\text{ for all }g\in H\}$ be the 
fixed subspace of $H$, and
$$
\cS=\{\fix(H)\,|\,H\text{ a subgroup of }G\}.
$$
Then $\fix(H)g=\fix(g^{-1}Hg)$, $\fix(H_1)\cap\fix(H_2)=\fix\langle
H_1,H_2\rangle$, where $\langle
H_1,H_2\rangle$ is the subgroup generated by the $H_i$,
and $V$ is the fixed space of the trivial subgroup. We will return to
this example in \S\ref{section:arrangementmonoids}.

\begin{definition}
\label{section:reflectionmonoids:definition100}
Let $G\subset GL(V)$ be a group and $\cS$ a system in $V$ for $G$. The
\emph{monoid of partial linear isomorphisms given by $G$ and $\cS$\/} 
is the submonoid of $\ml(V)$ defined by
$$
M(G,\cS):=\{g_X^{}\,|\, g\in G, X\in\cS\}.
$$
If $G$ is a reflection group then $M(G,\cS)$  is called a \emph{reflection monoid\/}.
\end{definition}

We will remind the reader of the definition of reflection group
in \S\ref{section:reflectiongroups}, and properly justify the terminology
``reflection monoid'' in \S\ref{section:generalities}.
Observe that the monoid structure on $M(G,\cS)$ is guaranteed
by (\ref{section:reflectionmonoids:equation100}) and (S1)-(S3).
If $\cS$ a system for $G$, $X\in\cS$, and $\ve:V\rightarrow V$ is the identity 
isomorphism, then $\ve_X^{}\in M(G,\cS)$, and thus every $X\in\cS$ is the
domain of some element of $M(G,\cS)$. Conversely, by (S1)-(S3) and 
(\ref{section:reflectionmonoids:equation100}), every element of $M(G,\cS)$ has domain
some element of $\cS$, so that $\cS$ is precisely the set of domains
of the partial isomorphisms in $M(G,\cS)$. If
$g_X^{}\in M(G,\cS)$ then $(g_X^{})^{-1}=g^{-1}_{Xg}\in M(G,\cS)$, and we have an inverse monoid
with units the $g\in G$ and idempotents the partial identities $\ve_X^{}$ for
$X\in\cS$. Moreover any $g_X^{}=\ve_X^{} g$, so $M(G,\cS)$ is factorizable.

Even with these very modest preliminaries, it is possible to prove a result with
non-trivial consequences:

\begin{theorem}\label{section:reflectionmonoids:result100}
Let $G\subset GL(V)$ be a finite group, $\cS$ a finite system
in $V$ for $G$, and $M(G,\cS)$ the resulting monoid of partial
linear isomorphisms. Then 
$$
|M(G,\cS)|=\sum_{X\in\cS}[G:G_X],
$$
where $G_X$ is the isotropy group of $X\in\cS$.
\end{theorem}

\begin{proof}
For $X\subset V$ let $M(X)$ be the set of $\aa\in M(G,\cS)$ with
domain $X$. Then $M(G,\cS)$ is the disjoint union of the 
$M(X)$, and as $\cS$ is precisely the set of domains of the $\aa\in M(G,\cS)$,
we have 
$|M(G,\cS)|=\sum_{X\in\cS}|M(X)|$. 
The 
elements of 
$M(X)$ are the partial isomorphisms obtained by restricting the elements
of $G$ to $X$, with $g_1,g_2\in G$ the same partial isomorphism if and only if they
lie in the same coset of the isotropy subgroup $G_X\subset G$. Thus
$|M(X)|$ is the index $[G:G_X]$ and the result follows.
\qed
\end{proof}

If $X,Y\in\cS$ lie in the same orbit of the $G$-action (S2) on $\cS$, then their
isotropy groups $G_X,G_Y$ are conjugate, and 
the sum in Theorem \ref{section:reflectionmonoids:result100} becomes 
\begin{equation}\label{section:reflectionmonoids:equation200}
|M(G,\cS)|=|G|\sum_{X\in\,\Om} \frac{n_X^{} }{|G_X|},
\end{equation}
where $\Om$ is a set of orbit representatives,
and $n_X^{}$ is the number of subspaces in the orbit containing $X$.
Most of our applications of Theorem \ref{section:reflectionmonoids:result100} will use the form
(\ref{section:reflectionmonoids:equation200}).

We end the section by recalling a result from semigroup theory. 
At several points in the paper we will want to identify a monoid of partial isomorphisms
with some pre-existing monoid in the literature. 
This is possible if the group of units are the same, the idempotents
are the same and the actions of the groups on the idempotents are the same:

\begin{proposition}\label{section:reflectionmonoids:result200}
Let $M=EG$ and $N=FH$  be factorizable inverse monoids,
and $\theta:G\rightarrow H$, $\varphi:E\rightarrow F$ homomorphisms
such that
\begin{itemize}
\item $\varphi$ is equivariant: $(geg^{-1})\varphi=(g\theta)(e\varphi)(g\theta)^{-1}$ for all
$g\in G$ and $e\in E$, and 
\item $\theta$ respects stablizers: $G_e\theta\subset H_{e\varphi}$ for all $e\in E$.
\end{itemize}
Then the map $\chi:M\rightarrow N$ given by
$(eg)\chi=(e\varphi)(g\theta)$ is a homomorphism.
Moreover,
$\chi$ is surjective if and only if $\theta,\varphi$ are surjective, and
$\chi$ is an isomorphism if and only if $\theta,\varphi$ are isomorphisms
with $G_e\theta=H_{e\varphi}$ for all $e\in E$.
\end{proposition}

We mention that this result also occurs in a preprint of D. Fitzgerald; as there, we leave
the straightforward proof to the reader.



\section{Two examples for the symmetric group $\SS_n$}
\label{section:symmetricgroup}

The representation of the 
symmetric group $\SS_n$ by permutation matrices leads to two interesting examples
of monoids of partial isomorphisms--both of which turn out to be reflection
monoids, and both of which can be identified with familiar monoids of partial permutations.

\subsection{The Boolean monoids}
\label{section:symmetricgroup:boolean}

Let $V$ be a Euclidean space
with basis $\{\xx_1,\ldots,\xx_n\}$ and inner product $(\xx_i,\xx_j)=\delta_{ij}$,
the Kronecker delta.
Let the symmetric group act on $V$ via $\xx_i\pi=\xx_{i\pi}$ for 
$\pi\in\SS_n$; we will abuse notation and write $\SS_n\subset GL(V)$,
identifying $\SS_n$ with its image under this representation.

For $J\subset I=\{1,\ldots,n\}$, let
\begin{equation}
  \label{section:symmetricgroup:equation25}
X(J)=\bigoplus_{j\in J}\R\xx_j\subset V,  
\end{equation}
and let $\BB$ be the collection of all such subspaces as $J$ ranges
over the subsets of $I$, with 
$X(\varnothing)=\mathbf{0}$. We have $X(I)=V$, 
$X(J)\pi=X(J\pi)$ for all $\pi\in\SS_n$, and
\begin{equation}
  \label{section:symmetricgroup:equation50}
X(J_1)\cap X(J_2)=X(J_1\cap J_2).  
\end{equation}
Indeed, partially ordering $\BB$ by inclusion, the map
$X(J)\rightarrow J$ is a lattice isomorphism $\BB\rightarrow\mathbf{2}^n$ to the 
Boolean lattice $\mathbf{2}^n$ of all subsets of $I$. 

The result is that $\BB$ is a system in $V$ for $\SS_n$, which in
honour of the lattice isomorphism above we will call the \emph{Boolean system\/}
for $\SS_n$. We form the associated monoid $M(\SS_n,\BB)$ and call it
the \emph{Boolean monoid\/}.

Clearly $X(J)$ has isotropy group $\SS_{I\setminus J}$.
Moreover, $\SS_I$ acts transitively on the subsets $J$ of a fixed size $k$,
so that 
a set $\Om$ of orbit representatives on $\BB$ is given 
by the $X(1,\ldots,k)$, with each orbit having size the number of 
$k$ element subsets of $I$. 
Plugging all of this into Theorem
\ref{section:reflectionmonoids:result100} and its alternative version
(\ref{section:reflectionmonoids:equation200}),
gives
\begin{equation}
  \label{section:symmetricgroup:equation100}
|M(\SS_n,\BB)|
=\sum_{J\subset I}[\SS_I:\SS_{I\setminus J}]
=|\SS_n|\sum_{k=0}^{n}{n\choose k}\frac{1}{|\SS_{n-k}|}
=\sum_{k=0}^{n}{n\choose k}^2 k!
\end{equation}
Readers from semigroup theory will recognize the formula on the right as the order
of the symmetric inverse monoid $\II_n$; readers from 
reflection groups will recognize the action of $\SS_n$ that led to it.

Neither is a coincidence: 
Let $\theta=\id:\SS_n\rightarrow\SS_n$ and $\varphi:\ve_{X(J)}\mapsto\ve_J^{}$
the isomorphism on the idempotents 
induced by the lattice isomorphism $\BB\rightarrow\mathbf{2}^n$---that
$\varphi$ is a homomorphism follows, for example, by 
(\ref{section:symmetricgroup:equation50}).
As $X(J)\pi=X(J\pi)$, we have $\pi^{-1}\ve_{X(J)}^{}\pi=\ve_{X(J\pi)}^{}$, and so 
$$
(\pi^{-1}\ve_{X(J)}^{}\pi)\varphi=\ve_{J\pi}^{}
=\pi^{-1}\ve_J^{}\pi=(\pi\theta)^{-1}(\ve_{X(J)}^{}\varphi)(\pi\theta),
$$ 
giving the
equivariance of $\varphi$. The stabilizer of the idempotent $e=\ve_{X(J)}^{}$ 
consists of those $\pi\in\SS_n$ such that $\xx_j\pi=\xx_j$ for all $j\in J$,
whereas for $e\varphi=\ve_J^{}$ we require $j\pi=j$. Proposition 
\ref{section:reflectionmonoids:result200} thus gives,

\begin{proposition}
\label{section:symmetricgroup:result300}
The map $\pi_{X(J)}^{}\mapsto \pi_J^{}$ is an isomorphism
$M(\SS_n,\BB)\rightarrow\II_n$ from the Boolean monoid to the symmetric inverse monoid. 
\end{proposition}

%

\subsection{The Coxeter arrangement monoids}
\label{section:symmetricgroup:coxeter}

We keep the same $\SS_n\subset GL(V)$ and notation from
\S\ref{section:symmetricgroup:boolean}, but switch 
to a more interesting system. For $1\leq i\not= j\leq n$,
let $\AA$ be the collection
of hyperplanes $H_{ij}=(\xx_i-\xx_j)^\perp\in V$, and $\HH=L(\AA)$
be the set of all possible intersections of elements
of $\AA$, with the null intersection taken to be $V$. This time we order
$\HH$ by \emph{reverse\/} inclusion, via which it
acquires the structure of a lattice, with the join of any two subspaces their intersection
and meet, the subspace generated by them.

Just as with the Boolean system $\BB$, 
we can identify $\HH$ with a well known combinatorial 
lattice. Recall that a partition of $I=\{1,\ldots,n\}$ is a collection
$\Lambda=\{\Lambda_1,\ldots,\Lambda_p\}$ of nonempty pairwise disjoint
subsets $\Lambda_i\subset I$, 
or \emph{blocks\/}, whose union is $I$. If $\lambda_i=|\Lambda_i|$ then 
$\lambda=\|\Lambda\|=(\lambda_1,\ldots,\lambda_p)$ is a partition of $n$, ie: 
the integers satisfy $\lambda_i\geq 1$ and $\sum\lambda_i=n$. 
Order the set $\Pi(n)$ of partitions of $I$ by
refinement: $\Lambda\leq\Lambda'$ if and only if every block of $\Lambda$
is contained in some block of $\Lambda'$.


The result is the partition lattice. 
It is not hard to show (see eg: \cite{Orlik92}*{Proposition 2.9})
that the map that sends 
the hyperplane $H_{ij}$ to
the partition with a single non-trivial
block $\Lambda_1=\{i,j\}$,
extends to 
a lattice isomorphism $\HH\rightarrow\Pi(n)$.
Indeed, if $X(\Lambda)\in\HH$ is the subspace,
\begin{equation}\label{section:symmetricgroup:equation200}
X(\Lambda)=\bigcap_{\lambda_k>1}\kern4pt\bigcap_{i,j\in\Lambda_k}H_{ij},
\end{equation}
then $X(\Lambda)\rightarrow\Lambda$ is the isomorphism. In particular
$\sum \aa_i\xx_i\in X(\Lambda)$ if and only if $\aa_i=\aa_j$ when $i,j$ lie in the
same block.
There is a faithful $\SS_n$-action on $\Pi(n)$ given by 
$\Lambda\pi=\{\Lambda_1\pi,\ldots,\Lambda_p\pi\}$, while the 
$\SS_n$-action on $\HH$ is given by $X(\Lambda)\pi=X(\Lambda\pi)$. 

As a consequence of this, and the fact that it is by definition closed
under intersection, we have $\HH$ is a system in $V$ for the symmetric group.
By virtue of its description as the intersection lattice 
for the ``arrangement'' $\AA$ of hyperplanes
$H_{ij}$, we call $\HH$ the \emph{Coxeter arrangement system\/} for $\SS_n$. We will
properly remind the reader about hyperplane arrangements in \S\ref{section:arrangementmonoids}. 
We remark that
the Boolean system $\BB$ of \S\ref{section:symmetricgroup:boolean} is also
an arrangement system, with $\BB=L(\AA)$, where $\AA$ consists
of the coordinate hyperplanes $\xx_i^\perp$. For reasons that will be made clearer in 
\S\ref{section:arrangementmonoids}, the $H_{ij}$ are a more natural collection of hyperplanes
to associate with the symmetric group than the $\xx_i^\perp$, so we
will reserve the arrangement terminology for this case.

We now apply Theorem \ref{section:reflectionmonoids:result100}. 
By (\ref{section:symmetricgroup:equation200}) and
the comments following it, we have $\xx\pi=\xx$ for all $\xx\in X(\Lambda)$ 
if and only if $\Lambda_i\pi=\Lambda_i$ for 
all $i$.
The isotropy group of the subspace $X(\Lambda)$ 
is thus isomorphic to a product of symmetric groups
$\SS_{\Lambda_1}\times\cdots\times\SS_{\Lambda_p}$, 
called a \emph{Young subgroup\/} of $\SS_n$, and Theorem 
\ref{section:reflectionmonoids:result100} becomes a sum, over all
partitions, of the indices of the resulting Young subgroups.

We can also be quite explicit:
for a partition $\Lambda$, let $b_i>0$ be the number of
$\lambda_j$ equal to $i$, and 
\begin{equation}
  \label{section:symmetricgroup:equation300}
b_\lambda=b_1!\ldots b_n!(1!)^{b_1}\ldots(n!)^{b_n}=b_1!\ldots b_n!
\lambda_1!\ldots\lambda_p!
\end{equation}

\begin{proposition}[\cite{Orlik92}*{Proposition 6.72}]
\label{section:symmetricgroup:result400}
In the action of the symmetric group $\SS_n$ on $\HH$,
two subspaces $X(\Lambda)$ and $X(\Lambda')$ lie in the same orbit if
and only if 
$\|\Lambda\|=\|\Lambda'\|$. The cardinality of the orbit of
the subspace $X(\Lambda)$ is $n!/b_\lambda$.
\end{proposition}

Plugging everything into version
(\ref{section:reflectionmonoids:equation200}) of Theorem
\ref{section:reflectionmonoids:result100}, including a summary of the 
discussion above, gives

\begin{theorem}\label{section:symmetricgroup:result500}
Let $\SS_n\subset GL(V)$ and $\HH=L(\AA)$ the intersection lattice
of the hyperplanes $H_{ij}$. 
Then
the Coxeter arrangement monoid $M(\SS_n,\HH)$ has order,
$$
|M(\SS_n,\HH)|
=\sum_\Lambda[\SS_I:\SS_{\Lambda_1}\times\cdots\times\SS_{\Lambda_p}]
=(n!)^2\sum_{\lambda}\frac{1}{b_\lambda \lambda_1!\ldots\lambda_p!},
$$ 
the first sum over all partitions $\Lambda$ of $I$, and the second
over all partitions $\lambda=(\lambda_1,\ldots,\lambda_p)$ of $n$, with
$b_\lambda$ given by (\ref{section:symmetricgroup:equation300}).
\end{theorem}

Theorem \ref{section:arrangementmonoids:result250} of 
\S\ref{section:arrangementmonoids:generalities} will
generalise this
result, replacing $\SS_n$ by an arbitrary finite reflection group.

The formula on the right hand side of Theorem \ref{section:symmetricgroup:result500}
may also ring a bell with the cognoscenti.
A \emph{uniform block permutation\/} 
is a bijection $\pi:\Lambda\rightarrow\Gamma$ between two partitions
of $I$. Thus, it is a bijection
$\pi:I\rightarrow I$, where the image of each block of $\Lambda$
is a block of $\Gamma$. If $\Lambda=\{\Lambda_1,\ldots,\Lambda_p\}$,
then up to a 
rearrangement of the blocks, $\Gamma=\{\Lambda_1\pi,\ldots,\Lambda_p\pi\}$
and we write $\lf\pi\rf_\Lambda^{}$ for this uniform block permutation,
noting that $\lf\pi\rf_\Lambda^{}=\lf\tau\rf_\Delta^{}$
if and only if $\Lambda=\Delta$ and $\Lambda_i\pi=\Delta_i\tau$ for all $i$.
We define an associative
product
$\lf\pi\rf_\Lambda^{}\lf\tau\rf_\Gamma^{}=\lf\pi\tau\rf_\Delta^{}$,
where $\Delta=\Lambda\vee\Gamma\pi^{-1}$ and 
$\vee$ is the join in the partition lattice (compare this expression
with the domain on the right hand side of
(\ref{section:reflectionmonoids:equation100})). This turns out to be a
factorizable 
inverse monoid, the \emph{monoid of uniform block permutations\/}
$\PP_n$ (see \cite{Aguiar08,Fitzgerald03,Kosuda00}).
Its group of units is $\SS_n$ and the idempotents are the 
$\lf\ve\rf_\Lambda^{}$ where $\ve:I\rightarrow I$ is the identity permutation.
We have $\pi^{-1}\lf\ve\rf_\Lambda^{}\pi=\lf\ve\rf_{\Lambda\pi}^{}$.


Let $\theta=\id:\SS_n\rightarrow\SS_n$ and $\varphi:\ve_{X(\Lambda)}^{}\mapsto\lf\ve\rf_\Lambda^{}$
the isomorphism on the idempotents 
induced by the lattice isomorphism $\HH\rightarrow\Pi(n)$:
that $\varphi$ is a homomorphism follows for example from
$X(\Lambda_1)\cap X(\Lambda_2)=X(\Lambda_1\vee\Lambda_2)$; remember that
$\HH$ is ordered by \emph{reverse\/} inclusion. Equivariance follows
from $X(\Lambda)\pi=X(\Lambda\pi)$ much as in the discussion preceding
Proposition \ref{section:symmetricgroup:result300}, as does the condition
on the idempotent stablizers. Thus,

\begin{proposition}\label{section:symmetricgroup:result600}
The map $\pi_{X(\Lambda)}^{}\mapsto\lf\pi\rf_\Lambda^{}$ is an isomorphism
$M(\SS_n,\HH)\rightarrow\PP_n$ from the Coxeter arrangement monoid to the
monoid of uniform block permutations.
\end{proposition}

\section{Reflection groups}
\label{section:reflectiongroups}

In this section we fix notation and leave the reader unfamiliar with
reflection groups to consult one of the standard references
\cites{Bourbaki02,Humphreys90,Kane01}. We have
generally followed \cite{Humphreys90}. 
For concreteness (as much as anything else) all the reflection groups
in this paper will be \emph{finite\/} and \emph{real\/}, that is,
subgroups $W\subset GL(V)$ generated by linear reflections of a real
vector space $V$.

Any (finite real) reflection group has the form $W(\Phi)=\langle
s_\vv | \vv\in\Phi\rangle$, where $\Phi\subset V$ is a root system and
$s_\vv$ the reflection in the hyperplane orthogonal to $\vv$.
%
The finite real reflection groups are, up to 
isomorphism, direct products of $W(\Phi)$ for $\Phi$ from a well known list of irreducible
root systems. 
These $\Phi$ fall into five infinite
families of types $A_{n-1},B_n,C_n,D_n$ (the classical systems) and $I_2(m)$,
and six exceptional cases of types $H_3, H_4, F_4, E_6,E_7$ and $E_8$.
Notable amongst these are the $\Phi$ whose
associated group $W(\Phi)$ is a finite crystallographic reflection, or
\emph{Weyl\/} group: these are the $W(\Phi)\subset GL(V)$ that leave invariant some
$\Z$-lattice $L\subset V$. The $W(\Phi)$ for $\Phi$ of type $I_2(m)$
are just the dihedral groups.

\begin{table}
\centering
\begin{tabular}{ccc}
\hline
Type of $\Phi$&Order of $W(\Phi)$&Root system $\Phi$\\\hline
\\
$A_{n-1}\,(n\geq 2)$&
$n!$&
$\{\xx_i-\xx_j \,\,(1\leq i\not= j\leq n)\}$\\
\\
$B_n\,(n\geq 2)$&
$2^n n!$&
$\{\pm\xx_i\,\,(1\leq i\leq n), \pm\xx_i\pm\xx_j\,\,(1\leq i<j\leq n)\}$\\
\\
$D_n\,(n\geq 4)$&
$2^{n-1}n!$&
$\{\pm\xx_i\pm\xx_j\,\,(1\leq i<j\leq n)\}$\\
\\
\hline
\end{tabular}\caption{Standard root systems 
$\Phi\subset V$
for the classical Weyl groups
\cite{Humphreys90}*{\S 2.10}.}\label{section:reflectiongroups:table1}
\end{table}

Table \ref{section:reflectiongroups:table1}
gives standard $\Phi$ for the classical Weyl
groups with 
$\{\xx_1,\ldots,$ $\xx_n\}$
an orthonormal basis for $V$.
The root systems of types $B$ and $C$ have the same symmetry, 
but different lengths of roots.
The associated Weyl groups are thus identical, and as it is these that 
ultimately concern us, we have given just the type $B$ system in 
Table \ref{section:reflectiongroups:table1}
(type $C$ has roots $\pm 2\xx_i$ rather than the $\pm\xx_i$). 
For convenience in expressing some of the formulae of \S\ref{section:booleanmonoids},
we extend the notation by adopting the additional  
conventions $A_{-1}=A_0=\varnothing$, $B_0=\varnothing,
B_1=\{\pm\xx_1\}$, and $D_0=D_1=\varnothing$,
$D_n=\{\pm\xx_i\pm\xx_j\,\,(1\leq i<j\leq n)\}$ for $n=2,3$.
Table \ref{section:reflectiongroups:table2} gives the orders of the
exceptional groups (where $W(G_2)$ is the group of type $I_2(6)$). We
will have no need for their root systems in this paper.

This paper started with the monoids $M(G,\cS)$ for $G$ an arbitrary
(linear) group. 
One of the reasons to focus on the case that $G$ is a reflection group
is that the isotropy groups
$G_X$ are very often also reflection groups, and this makes the calculation in 
(\ref{section:reflectionmonoids:equation200}) do-able. 
A theorem of Steinberg \cite{Steinberg60}*{Theorem 1.5}
asserts that for $G=W(\Phi)$ and $X\subset V$ any subspace, the
isotropy group
$W(\Phi)_X$ is generated by the reflections $s_\vv$ for
$\vv\in\Phi\cap X^\perp$.

\section{The Boolean reflection monoids}
\label{section:booleanmonoids}

In \S\ref{section:symmetricgroup:boolean} we considered the
permutation action of $\SS_n$ on a Euclidean space $V$ of dimension
$n$: we now know that this is nothing more than a well known realization of $\SS_n$ as the
reflection group $W(A_{n-1})$. Indeed
the map $s_{\xx_i-\xx_j}\mapsto (i,j)$ induces an isomorphism
$W(A_{n-1})\rightarrow\SS_n$, which we write as $g(\pi)\mapsto\pi$.
Moreover, the $W(A_{n-1})$-action on the subspaces $X(J)\in\BB$ is just $X(J)g(\pi)=X(J\pi)$.

The Boolean monoid $M(\SS_n,\BB)$ of
\S\ref{section:symmetricgroup:boolean} is thus a reflection monoid,
which we will denote as $M(A_{n-1},\BB)$ from now on, and the isomorphism
$W(A_{n-1})\cong\SS_n$ extends to an isomorphism $M(A_{n-1},\BB)\cong\II_n$.

The moral of this section is that the Boolean $\BB$ of
\S\ref{section:symmetricgroup:boolean} is a system for 
all the classical Weyl groups of types $A,B,D$, and the orders of
the resulting reflection monoids 
can be determined in a nice uniform way.
Moreover, these Weyl groups
have well known alternative descriptions as certain groups of
permutations, and so too the 
resulting reflection monoids,
at least in types $A$ and $B$, have descriptions as naturally occurring monoids
of permutations.

\begin{table}
\centering
\begin{tabular}{cc}
\hline
Type of $\Phi$&Order of $W(\Phi)$\\\hline
$G_2$\vrule width 0 mm height 5 mm depth 0 mm&
$2^3\,3$
\\
$F_4$&
$2^7\,3^2$
\\
$E_6$&
$2^7\,3^4\,5$
\\
$E_7$&
$2^{10}\,3^4\,5\,7$
\\
$E_8$&
$2^{14}\,3^5\,5^2\,7$\vrule width 0 mm height 0 mm depth 3 mm\\
\hline
\end{tabular}\caption{Exceptional Weyl groups and their
  orders.}\label{section:reflectiongroups:table2} 
\end{table}

First, the alternative descriptions of the reflection groups in 
types $B$ and $D$. As usual $I$ is a set and $-I=\{-x\,|\,x\in I\}$ 
a set with the same cardinality. The group
$\SS_{\pm I}$ of {\em signed permutations\/}
of $I$
is $\SS_{\pm I}=\{\pi\in\SS_{I\cup-I}\,|\,(-x)\pi=-(x\pi)\}$, where
$\SS_{\pm n}$ has the obvious meaning.
For $x\in I$, let $|x|=\{x,-x\}$, and $|I|=\{|x|:x\in I\}$.
If $\pi\in\SS_{\pm I}$, define $|\pi|\in\SS_{|I|}$ by
$$
|x||\pi|=|y|\Leftrightarrow\{x\pi,-x\pi\}=\{y,-y\}.
$$
Then the map $\pi\rightarrow|\pi|$ is a surjective homomorphism 
$|\cdot|:\SS_{\pm I}\rightarrow\SS_{|I|}\cong\SS_I$. 

A signed permutation $\pi$ is \emph{even\/} if the number of $x\in I$ with
$x\pi\in-I$ is even, and the even signed permutations form a subgroup
$\SS_{\pm I}^e$ of index two in $\SS_{\pm I}$. Indeed, if $\tau_x$ is the signed
transposition $(x,-x)$, then $\{1,\tau_x\}$ are coset representatives 
for $\SS_{\pm I}^e$ in $\SS_{\pm I}$. In particular, as $|\tau_x|=1$, 
restriction gives a surjective homomorphism
$|\cdot|:\SS_{\pm I}^e\rightarrow\SS_I$.

There are isomorphisms $W(B_n)\rightarrow\SS_{\pm n}$ induced by
$$
s_{\xx_i-\xx_j}\mapsto(i,j)(-i,-j)\text{ and }s_{\xx_i}\mapsto (i,-i)
$$
and $W(D_n)\rightarrow\SS_{\pm n}^e$ induced by
$s_{\xx_i-\xx_j}\mapsto(i,j)(-i,-j)$ and $s_{\xx_i+\xx_j}\mapsto (i,-j)(-i,j)$.
As above, we write $g(\pi)$ for the element of $W(B_n)$ or $W(D_n)$ corresponding
to $\pi\in\SS_{\pm n}$ or $\SS_{\pm n}^e$.

Now let $V$ be Euclidean with orthonormal basis $\{\xx_1\ldots,\xx_n\}$, 
$I=\{1,\ldots,n\}$, and 
\begin{equation}
  \label{section:booleanmonoids:equation100}
\BB=\{X(J)\,|\,J\subset I\},  
\end{equation}
the subspaces 
from (\ref{section:symmetricgroup:equation25}). 
Using these descriptions of the reflection groups $W(\Phi)$ for $\Phi=A_{n-1},B_n$ and $D_n$,
we have
\begin{equation}
  \label{section:booleanmonoids:equation150}
X(J)g(\pi)=X(J|\pi|)  
\end{equation}
for $g(\pi)\in W(\Phi)$, as well as the $X(I)=V$ and $X(J_1)\cap X(J_2)=X(J_1\cap J_2)$ that
we had in \S\ref{section:symmetricgroup:boolean}. Thus $\BB$ is a system in $V$ for
$W(\Phi)$, which we continue to call the \emph{Boolean system}. We
write $M(\Phi,\BB)$ instead of $M(W(\Phi),\BB)$, and call these the
\emph{Boolean reflection monoids\/} of types $A_{n-1}, B_n$ or $D_n$. Note that 
$\BB$ is not a system for any of the exceptional $W(\Phi)$.

Now to their orders.
Let $\Phi=\Phi_n$ be a root system of type $A_{n-1}, B_n$ or $D_n$ as in Table 
\ref{section:reflectiongroups:table1}. For $J\subset I$, we write $\Phi_J$
for $\Phi\cap X(J)$. By Steinberg's theorem the isotropy group of the
subspace $X(J)$ is generated 
by the $s_\vv$ with $\vv\in\Phi\cap X(J)^\perp=\Phi_{I\setminus J}$.
As the homomorphism $|\cdot|$ maps $\SS_n,\SS_{\pm n}$ and $\SS_{\pm n}^e$ onto
$\SS_n$, the action (\ref{section:booleanmonoids:equation150}) of
$W(\Phi)$ on $\BB$ is transitive, for 
fixed $k$, on the $X(J)$ of dimension $k$. Thus:

\begin{theorem}\label{section:booleanmonoids:result100}
Let $\Phi_n$ be a root system of type $A_{n-1},B_n$ or $D_n$ as in Table
\ref{section:reflectiongroups:table1} and $\BB$ the Boolean system 
(\ref{section:booleanmonoids:equation100})
for $W(\Phi_n)$.
Then the Boolean reflection monoids have orders,
$$
|M(\Phi_n,\BB)|
=\sum_{J\subset I}[W(\Phi):W(\Phi_{I\setminus J})]
=|W(\Phi_n)|\sum_{k=0}^n 
{n\choose k}\frac{1}{|W(\Phi_{n-k})|}.
$$
\end{theorem}

Compare Theorem \ref{section:booleanmonoids:result100} with
(\ref{section:symmetricgroup:equation100}).
By the conventions of \S\ref{section:reflectiongroups}
we have $|W(A_k)|=(k+1)!$, $|W(B_k)|=2^k k!$, $|W(D_0)|=1$,
and $|W(D_k)|=2^{k-1} k!$ for $k>1$, giving the explicit versions,
$$
\begin{tabular}{c|ccc}
\hline
$\Phi_n$&$A_{n-1}$&$B_n$&$D_n$\\\hline
&&&\\
$|M(\Phi_n,\BB)|$
&\vrule width 2 mm height 0 mm depth 0 mm
${\ds \sum_{k=0}^{n}{n\choose k}^2
k!}$
&\vrule width 2 mm height 0 mm depth 0 mm
${\ds \sum_{k=0}^{n}2^k{n\choose k}^2
k!}$
&\vrule width 3 mm height 0 mm depth 0 mm
${\ds 2^{n-1}n!+\sum_{k=1}^{n}2^k{n\choose k}^2
k!}$\\
&&&\\\hline
\end{tabular}
$$

The dichotomy $W(A_{n-1})\cong\SS_n, M(A_{n-1},\BB)\cong\II_n$
has a type $B$ version, for which we need an inverse monoid to play the
role of $\SS_{\pm n}$.
This is the \emph{monoid of partial signed permutations\/} of $I$:
$$
\II_{\pm I}:=\{\pi\in\II_{I\cup -I}\,|\, (-x)\pi=-(x\pi)\text{ and }x\in\dom\,\pi
\Leftrightarrow -x\in\dom\,\pi\},
$$
with $\II_{\pm n}$ having the obvious meaning.
Every element of $\II_{\pm n}$ has the form $\pi_X^{}, (X=J\cup -J)$ for some
signed permutation $\pi$ and $J\subset I$. Thus $\II_{\pm n}$ has
units $\SS_{\pm I}$ and idempotents the $\ve_X^{}, (X=J\cup -J)$ with 
$\ve:I\rightarrow I$ the identity map.

Let $\theta:W(B_n)\rightarrow\SS_{\pm n}$ be the isomorphism $g(\pi)\mapsto\pi$ described above,
and 
$$\varphi:\ve_{X(J)}^{}\mapsto\ve_X^{}\,(X=J\cup -J)$$
the isomorphism on the idempotents induced
by the lattice isomorphism $\BB\rightarrow\mathbf{2}^n$. Observe that
if $\pi\in\SS_{\pm n}$ 
then $\pi^{-1}\ve_X^{}\pi=\ve_{X|\pi|}^{}, (X|\pi|=J|\pi|\cup -J|\pi|)$, and the equivariance of
$\varphi$ follows from this and $X(J)\pi=X(J|\pi|)$. Thus, another application of
Proposition \ref{section:reflectionmonoids:result200} gives,

\begin{proposition}\label{examples:permutationmonoidsresult200}
The map $g(\pi)_{X(J)}^{}\mapsto \pi_X^{}, (X=J\cup -J)$ is an isomorphism
$M(B_n,\BB)\rightarrow\II_{\pm n}$ from the Boolean monoid of type $B$ to the 
monoid of partial signed permutations.
\end{proposition}

Thus we have the pair $W(B_n)\cong\SS_{\pm n}$ and $M(B_n,\BB)\cong\II_{\pm n}$, to go
with the one in type $A$.
What about a pair $W(D_n)\cong\SS_{\pm n}^e$ and 
$M(D_n,\BB)\cong\II_{\pm n}^e$,
or some such? 
The problem is that one can show, by thinking in terms
of partial signed permutations, that the non-units of $M(B_n,\BB)$ and $M(D_n,\BB)$
are the same (which is why the orders of these reflection monoids are identical
except for the $k=0$ terms). This makes a nice interpretation of $M(D_n,\BB)$ in
terms of ``even signed permutations'' unlikely.

\section{The Coxeter arrangement monoids}
\label{section:arrangementmonoids}

Just as \S\ref{section:booleanmonoids} generalizes the Boolean monoid
$M(\SS_n,\BB)$ of \S\ref{section:symmetricgroup:boolean}, replacing
$\SS_n$ by a classical Weyl group, so now we generalize the Coxeter
arrangement monoid $M(\SS_n,\HH)$ of
\S\ref{section:symmetricgroup:coxeter}, replacing $\SS_n$ by an
arbitrary finite reflection group.

\subsection{Generalities}
\label{section:arrangementmonoids:generalities}

Steinberg's Theorem (\S\ref{section:reflectiongroups}) provides a good
reason to study reflection monoids, rather than just monoids of
partial isomorphisms. 
Another reason is that 
the system 
$$
\cS=\{\fix(H)\,|\,H\text{ a subgroup of }G\},
$$
of \S\ref{section:reflectionmonoids}
has a particularly nice combinatorial structure when $G$ is a
reflection group.
These systems and the resulting reflection monoids, especially when $G$ is a Weyl group,
are the subject of this section. Recall that $V$ is a finite dimensional
real space and $W\subset GL(V)$
a finite reflection group.

A \emph{hyperplane arrangement\/} $\AA$ is a finite collection of hyperplanes
in $V$. General references are \cites{Orlik92,Zaslavsky75},
where the hyperplanes can be affine, but we restrict ourselves to 
linear arrangements.
An important combinatorial invariant for $\AA$ is the
intersection lattice $L(\AA)$--the set of all possible intersections of elements
of $\AA$, ordered by reverse inclusion, and with the null intersection taken to be
the ambient space $V$. What results is a lattice \cite{Orlik92}*{\S 2.1},
with unique minimal element the space $V$. If the $\AA$ are the reflecting hyperplanes of a 
reflection group $W\subset GL(V)$, then we have a \emph{reflection\/}
or \emph{Coxeter arrangement\/}.

The intersection lattice $L(\AA)$ of a Coxeter arrangement is a system of subspaces
for the associated reflection group $W$:
if $X\in\AA$ and $s_X^{}\in W$ is the
reflection in $X$, then for $g\in W$ we have $s_{Xg}^{}=g^{-1}s_X^{}g$, 
and so $Xg\in\AA$. Thus $\AA\,W=\AA$, extending to $L(\AA)\,W=L(\AA)$.
We will write $\HH$ for $L(\AA)$, calling it the \emph{reflection\/}
or \emph{Coxeter arrangement system\/}, and reserving $\BB$ for the Boolean system.
We call the resulting $M(W,\HH)$ the \emph{reflection\/}
or \emph{Coxeter arrangement monoid\/} of $\AA$ (or $W$).
If $W=W(\Phi)$, we write $\HH(\Phi)$ for $\HH$, 
and $M(\Phi,\HH)$ for $M(W(\Phi),\HH)$, 
observing that $M(A_{n-1},\HH)$ is the
monoid $M(\SS_n,\HH)$, or the monoid of uniform block
permutations, of \S\ref{section:symmetricgroup:coxeter}.

\begin{lemma}
Let $W\subset GL(V)$ be a finite reflection group with reflecting hyperplanes $\AA$
and Coxeter arrangement system $\HH=L(\AA)$.
Then $\HH=\{\fix(H)\,|\,H\text{ a subgroup of }W\}$.
\end{lemma}

\begin{proof}
Write $\FF$ for the system of fixed subspaces.
It is not hard to show using Steinberg's 
theorem and induction on the order of $W$ (see, eg: \cite{Orlik92}*{Theorem 6.27}), that
$\text{Fix}(g)\in\HH$ for all $g\in W$, where $\text{Fix}(g)$ is the fixed
subspace of the element $g$. As 
$\text{Fix}(H)=\bigcap_{H}\text{Fix}(h)$ and $H$ is finite, we
get $\FF\subset\HH$. 
Moreover, if $s$ is a reflection  then $\text{Fix}(s)$ is the reflecting 
hyperplane of $s$, so that $\AA\subset\FF\subset\HH$. As every element of $\HH$
is an intersection of elements of $\AA$, and $\FF$ is closed under $\bigcap$,
we have $\FF=\HH$.
\qed
\end{proof}

We will see in \S\ref{section:generalities} that a monoid isomorphism
$M(G,\cS)\rightarrow M(G',\cS')$ 
induces a poset isomorphism $\cS\rightarrow\cS'$, with the subspaces ordered by inclusion.
A comparison of the number of $k$-dimensional subspaces in $\HH$ and
$\BB$ shows that for a given $W$, there can be no isomorphism between
the Boolean and Coxeter arrangement monoids (see
\cite{Orlik92}*{\S6.4} for the number of subspaces in $\HH$ given in
terms of Stirling numbers of the second kind).

Now to orders: let $W=W(\Phi)$ and
$\Delta\subset\Phi$ a simple system. If $I\subset\Delta$, let
$W_I=\langle s_{\xx} \,|\,\xx\in I\rangle$ be the resulting 
\emph{special parabolic\/} subgroup, with
a \emph{parabolic\/} subgroup being any $W$-conjugate of a special
parabolic (see \cite{Humphreys90}*{\S1.10}). The parabolic subgroups
are thus
parametrised by the pairs $I,w$ with $I\subset\Delta$ and $w$ a
(right) coset representative for $W_I$ in $W$.

The parabolics in $\SS_I$ are just the Young subgroups
$\SS_{\Lambda_1}\times\cdots\times\SS_{\Lambda_p}$ for
$\Lambda=\{\Lambda_1,\ldots,\Lambda_p\}$ a partition of $I$. Theorem
\ref{section:symmetricgroup:result500} gives the order of
$M(\SS_n,\HH)$ as the sum of the indices of these, and indeed this is
the case in general:

\begin{theorem}\label{section:arrangementmonoids:result250}
Let $W\subset GL(V)$ be a finite reflection group 
with Coxeter arrangement system
$\HH$. Then the Coxeter arrangement monoid $M(W,\HH)$ has order
the sum of the indices of the parabolic subgroups 
of $W$.
\end{theorem}

\begin{proof}
By \cite{Kane01}*{Theorem 5.2} the isotropy groups $W_X$ are
parabolic, so it suffices to show that every parabolic
arises as an isotropy group $W_X$ for some $X\in\HH$, and that distinct subspaces
in $\HH$ have distinct isotropy groups.

The second of these requires only elementary arguments: if $X,Y$ are any subspaces of $V$
with $W_X=W_Y=W_\aa$, then $W_\aa$ also fixes $X+Y$ pointwise. If $X$ and $Y$ are distinct,
with one not contained in the other, then $X$ say, is a proper subspace of $X+Y$.
Thus it suffices to show that $X\subsetneqq Y$ have distinct isotropy
groups for $X,Y\in\HH$. Suppose otherwise, and recalling that
$X,Y\in\HH$ are intersections of reflecting hyperplanes of $W$, write
$X=Y\cap H_1\cap\cdots\cap H_k$ with the $H_i$ reflecting hyperplanes of $W$
and $Y\not\subset H_i$ for any $i$. In particular $W_{H_i}=\langle s_i\rangle$
with $s_i$ the reflection in $H_i$, and $s_i\not\in W_Y$
\cite{Humphreys90}*{Theorem 1.10}. But we also have
$\langle W_Y,s_1,\ldots,s_k\rangle\subset W_X=W_Y$, so that the $s_i\in W_Y$,
a contradiction.

Next we reduce to the case where the roots $\Phi$ span $V$.
Let $U=\sum_\Phi\R\xx\subset V$ and decompose $V=U\oplus U^\perp$. 
Then $U$ is a $W$-invariant subspace of $V$ with 
the complement $U^\perp$ fixed pointwise,
and $W$ is a finite reflection group in $GL(U)$ with reflecting hyperplanes the
$\xx^\perp\cap U$ for $\xx\in\Phi$, and associated Coxeter arrangement
system $\HH_U:=\{X\cap U\,|\,X\in\HH\}$. The map $X\mapsto X\cap U$ is a lattice
isomorphism $\HH\rightarrow\HH_U$ which induces an isomorphism between
the idempotents of $M(W,\HH)$ and $M(W,\HH_U)$. An application of Proposition 
\ref{section:reflectionmonoids:result200} 
then gives an isomorphism between these two reflection monoids.
Henceforth then, we will assume that $\sum_\Phi\R\xx=V$. 

We now
remind the reader of the \emph{Coxeter complex\/} in $V$, whose cells
have isotropy groups easily identified with specific parabolics
(see \cite{Humphreys90}*{1.15}). 
If $I\subset\Delta$ let $C_I$ be the subset of $V$ given by:
$$
C_I:=\{\vv\in V\,|\,(\xx,\vv)=0\text{ for }\xx\in I\text{ and }
(\xx,\vv)>0\text{ for }\xx\in\Delta\setminus I\}
=\bigcap_{\xx\in I}\xx^\perp\cap\bigcap_{\xx\in\Delta\setminus I}\xx^{>0},
$$
where $\xx^{>0}$ is the open half space consisting of those $\vv$
with $(\xx,\vv)>0$. The Coxeter complex 
$\Sigma$ has codimension $k$ cells
the subsets $C_Iw\subset V$ where $w\in W$ and $|I|=k$, and
the isotropy group of
the cell $C_Iw$ is the parabolic $w^{-1}W_I w$. 
If $\Delta'\subset\Phi$ is some other simple system then 
there
is a $w\in W$ with $\Delta'=\Delta w$. Hence if $I'\subset\Delta'$, then
\begin{equation}
  \label{section:arrangementmonoids:equation100}
C_{I'}=\bigcap_{\xx\in I'}\xx^\perp\cap\bigcap_{\xx\in\Delta'\setminus I'}\xx^{>0},  
\end{equation}
is the cell $C_Iw\in\Sigma$ for $I=I'w^{-1}$, and every cell of $\Sigma$
has this form. Thus the parabolic $w^{-1}W_I w$ is the isotropy group
$W_Y$ for the $Y=C_{I'}=C_Iw$ of (\ref{section:arrangementmonoids:equation100}).
If $X=\bigcap_{\xx\in I'}\xx^\perp$
then $X\in\HH$ and $Y\subset X$ is an open subset.
In particular, $Y$ spans $X$, so that $W_Y=W_X$. Hence every parabolic arises as a $W_X$ for
some $X\in\HH$.
\qed  
\end{proof}

\subsection{Coxeter arrangement monoids of type $B$}
\label{section:arrangementmonoids:typeB}

Just as in Theorem \ref{section:symmetricgroup:result500}, we can more explicit about
the orders of the Coxeter arrangement monoids. 
The material here and in \S\ref{section:arrangementmonoids:typeD} is adapted from
\cite{Orlik92}*{\S6.4}.

We build a combinatorial model for the 
Coxeter arrangement system $\HH(B_n)$, much as the partition lattice $\Pi(n)$ models
the system in type $A$. If $I=\{1,\ldots,n\}$, then a \emph{coupled
partition\/} $\Lambda$ of $I$ is a collection,
\begin{equation}
  \label{section:arrangementmonoids:equation200}
\Lambda=\{\Lambda_{11}+\Lambda_{12},\ldots,\Lambda_{q1}+\Lambda_{q2},
\Lambda_1,\ldots,\Lambda_p\},  
\end{equation}
of non-empty pairwise disjoint subsets whose union is $I$. The $\Lambda_{ij}$
and $\Lambda_i$ are \emph{blocks\/}, with $\Lambda_{i1}+\Lambda_{i2}$ a
\emph{coupled block\/}. The $+$ sign is purely formal, indicating
that these two blocks have been coupled. Thus, a coupled partition is just a partition
with some extra structure.
The coupled partition is completely determined by the blocks and the couplings, so that
reordering the blocks, the coupled blocks
or even the blocks within a coupled block, gives the same coupled partition.
If $\lambda_{ij}=|\Lambda_{ij}|$ then let $\lambda=\|\Lambda\|=(\lambda_{11}+\lambda_{12},
\ldots,\lambda_{q1}+\lambda_{q2},\lambda_1,\ldots,\lambda_p)$ be 
the resulting partition of $n$, where now $+$ really \emph{does\/} mean $+$.

Let $\TT$ be the set of pairs $(\Delta,\Lambda)$ where $\Delta\subset I$
and $\Lambda$ is a coupled partition of $I\setminus\Delta$. Define a
relation on $\TT$ by
$(\Delta,\Lambda)\leq(\Delta',\Lambda')$ if and only if 
\begin{itemize}
\item $\Delta\subset\Delta'$;
\item each (uncoupled) block $\Lambda_i\in\Lambda$ is either contained in $\Delta'$, or an
(uncoupled) block
$\Lambda'_j\in\Lambda'$ or a block $\Lambda'_{ij}$ of a couple $\in\Lambda'$;
\item each couple $\Lambda_{i1}+\Lambda_{i2}$ is either contained in $\Delta'$,
or $\Lambda_{i1}\subset\Lambda'_{j1}$ and $\Lambda_{i2}\subset\Lambda'_{j2}$
for some couple $\Lambda'_{j1}+\Lambda'_{j2}\in\Lambda'$.
\end{itemize}

Now let $V$ be Euclidean with orthonormal basis $\{\xx_1\ldots,\xx_n\}$.
For $(\Delta,\Lambda)\in\TT$, let $X(\Delta,\Lambda)\subset V$ be the subspace consisting
of those $\xx=\sum\aa_i\xx_i$ where $\aa_i=0$ for $i\in\Delta$, $\aa_i=\aa_j$ if 
$i,j$ lie in the same block of $\Lambda$ (either uncoupled or in a
couple) and $\aa_i=-\aa_j$ if $i,j$ lie in different 
blocks of the same coupled block. By the results of \cite{Orlik92}*{\S6.4}, the 
Coxeter arrangement system $\HH(B_n)$ has elements the
subspaces $X(\Delta,\Lambda)$ for $(\Delta,\Lambda)\in\TT$.

\begin{proposition}\label{section:arrangementmonoids:result300}
The map $f: (\Delta,\Lambda)\mapsto X(\Delta,\Lambda)$ is a bijection
from $\TT$ to $\HH(B_n)$ with
$(\Delta,\Lambda)\leq(\Delta',\Lambda')$ if and only if $X(\Delta',\Lambda')
\subset X(\Delta,\Lambda)$. In particular, $\TT$ with the relation $\leq$ 
defined above is a lattice and $f$
is a lattice isomorphism
$\TT\rightarrow\HH(B_n)$.
\end{proposition}

Orlik and Terao \cite{Orlik92}*{\S6.4} parametrize the subspaces in $\HH(B_n)$
using triples consisting of a subset $\Delta$, a \emph{partition\/}
$\Lambda$ of $I\setminus\Delta$ 
and a $\Gamma\subset I$, although these triples do not have a lattice structure.

\begin{proof}
That $f$ is a bijection is a straight notational translation of the 
results of \cite{Orlik92}*{\S6.4}. If $(\Delta,\Lambda)\leq(\Delta',\Lambda')$ and
$\xx=\sum\aa_i\xx_i\in X(\Delta',\Lambda')$, then it is easy to show that $\xx$ satisfies the
conditions for being an element of $X(\Delta,\Lambda)$: we have $\aa_i=0$ for
$i\in\Delta'$, so that $\aa_i=0$ when $i\in\Delta$, as $\Delta\subset\Delta'$,
and so on. Conversely, by looking at the coordinates of the elements of
$X(\Delta,\Lambda)$, it is easy to see that any subspace of the form
$X(\Delta',\Lambda')$ must satisfy $(\Delta,\Lambda)\leq(\Delta',\Lambda')$.
\qed
\end{proof}

We saw in \S\ref{section:booleanmonoids} that $W(B_n)\cong\SS_{\pm n}$.
There is another alternative (and well known) description of $W(B_n)$ that is useful 
in the current context.
Let $\mathbf{2}^I$ be the subsets of $I$, but now an Abelian group
under symmetric difference $S\vartriangle T:=(S\cup T)\setminus (S\cap T)$.
The symmetric group $\SS_I$ acts
on $\mathbf{2}^I$ via $T\mapsto T\pi$, ($\pi\in\SS_n$ and $T\subset I$),
and we form the 
semi-direct product $\SS_I\ltimes\mathbf{2}^I$, in which every element has a 
unique expression as a pair $(\pi,T)$ with $\pi\in\SS_I$, $T\subset I$. 
The map
$s_{\xx_i-\xx_j}\mapsto (i,j)\in\SS_n$, $s_{\xx_i}\mapsto\{i\}\in\mathbf{2}^n$ 
induces an isomorphism 
\begin{equation}
  \label{section:arrangementmonoids:equation300}
W(B_n)\rightarrow\SS_n\ltimes\mathbf{2}^n,
\end{equation}
and we write $g(\pi,T)\in W(B_n)$ for the element mapping to $(\pi,T)$.
If $J\subset I$, then $\SS_J\ltimes\mathbf{2}^J$ is naturally a subgroup
of $\SS_I\ltimes\mathbf{2}^I$.

For $T\in\mathbf{2}^I$ and $J\subset I$, let $J^+=J\cap T$ and 
$J^-=J\setminus T$, decomposing $J$ as a disjoint union
$J=J^+\cup J^-$. If $\Lambda$ is the coupled partition
(\ref{section:arrangementmonoids:equation200}) 
and $\pi\in\SS_n$, then let 
$$
\Lambda\pi:=\{\ldots,\Lambda_{i1}\pi+\Lambda_{i2}\pi,\ldots,\Lambda_i\pi,\ldots\},
$$
and if $T\in\mathbf{2}^I$, then let 
$$
\Lambda^T:=\{\ldots,(\Lambda_{i1}^-\cup\Lambda_{i2}^+)
+(\Lambda_{i2}^-\cup\Lambda_{i1}^+),\ldots,\Lambda_i^-+\Lambda_i^+,\ldots\},
$$
with the convention $\Lambda+\varnothing=\varnothing+\Lambda=\Lambda$.
Define an action of $\SS_I\ltimes\mathbf{2}^I$ on the lattice $\TT$ by
\begin{equation}
  \label{section:arrangementmonoids:equation400}
(\Delta,\Lambda)(\pi,T)=
(\Delta\pi,\Lambda^T\pi).
\end{equation}
Just as $\SS_n\ltimes\mathbf{2}^n$ models $W(B_n)$ and $\TT$ models $\HH(B_n)$, so
(\ref{section:arrangementmonoids:equation400}) models the action of $W(B_n)$ on $\HH(B_n)$:
we have $X(\Delta,\Lambda)g(\pi,T)=X(\Delta\pi,\Lambda^T\pi)$. In
particular, the 
lattice isomorphism of Proposition \ref{section:arrangementmonoids:result300}
is equivariant with respect to the $\SS_n\ltimes\mathbf{2}^n$ action on $\TT$ and the 
$W(B_n)$ action on $\HH(B_n)$.

This observation allows us to give a combinatorial version of the
Coxeter arrangement monoid $M(B_n,\HH)$. Its elements are ``uniform
block signed permutations'' of the elements of $\TT$ with the action
just described, and may be written in the form 
$\lf\pi,T\rf_{(\Delta,\Lambda)}$ where
$(\pi,T)\in\SS_I\ltimes\mathbf{2}^I$ and $(\Delta,\Lambda)\in\TT$. We
have
$\lf\pi,T\rf_{(\Delta,\Lambda)}=\lf\pi',T'\rf_{(\Delta',\Lambda')}$ if and
only if $\Delta=\Delta'$, $\Lambda=\Lambda'$, $\Delta\pi=\Delta'\pi'$
and $\Lambda_i\pi=\Lambda_i'\pi'$ for all $i$. The product is defined
by 
$$
\lf\pi,T\rf_{(\Delta,\Lambda)}\lf\pi',T'\rf_{(\Delta',\Lambda')}
=\lf(\pi,T)(\pi',T')\rf_{(\Gamma,\Upsilon)},
$$
where
$(\Gamma,\Upsilon)=(\Delta,\Lambda)\vee(\Delta',\Lambda')(\pi,T)^{-1}$,
with $\vee$ the join in the lattice $\TT$. We leave the diligent
reader to verify that this definition does indeed give a monoid
isomorphic to $M(B_n,\HH)$, and content ourselves with the observation
that this example illustrates the advantage of the geometric approach over the
combinatorial one.

For a coupled block $\Lambda_{1}+\Lambda_{2}\subset I$, let
$\BBB_{\Lambda_{1}+\Lambda_{2}}\subset
\SS_{\Lambda_{1}+\Lambda_{2}}\ltimes\mathbf{2}^{\Lambda_{1}+\Lambda_{2}}$  
be the subgroup consisting of those $(\pi,T)$ in
that leave each block of the couple invariant under the 
restriction of the action (\ref{section:arrangementmonoids:equation400}). 
Precisely, we require $(\Lambda_{1}^+\cup\Lambda_{2}^-)\pi=\Lambda_{1}$, from which it
follows that $(\Lambda_{2}^+\cup\Lambda_{1}^-)\pi=\Lambda_{2}$.

\begin{lemma}\label{section:arrangementmonoids:result350}
Let $\lambda_{i}=|\Lambda_{i}|$. Then the group $\BBB_{\Lambda_{1}+\Lambda_{2}}$ 
has order $\lambda_{1}!\lambda_{2}!c(\lambda_{1}\lambda_{2})$, where
$$
c(\lambda_{1},\lambda_{2})=\kern-3mm\sum_{j=0}^{\min\{\lambda_{1},\lambda_{2}\}}
{\lambda_{1}\choose j}{\lambda_{2}\choose j}.
$$
\end{lemma}

\begin{proof}
To have $(\Lambda_{1}^+\cup\Lambda_{2}^-)\pi=\Lambda_{1}$, it is clearly necessary that
$\Lambda_{1}^+\cup\Lambda_{2}^-$ and $\Lambda_1$ have the same cardinality, and 
on the other hand, if this is so, then $\pi$ can be a bijection 
extending any bijection
$\Lambda_{1}^+\cup\Lambda_{2}^-\rightarrow\Lambda_1$. The $T$ for 
which these two sets have the same size are precisely those 
with $|T\cap\Lambda_1|=|T\cap\Lambda_2|$,
of which there are
$c(\lambda_1,\lambda_2)$, and for each one there are $\lambda_1!$ bijections 
$\Lambda_{1}^+\cup\Lambda_{2}^-\rightarrow\Lambda_1$, each one in turn extendable to
$\lambda_2!$ bijections $\pi:\Lambda_{1}+\Lambda_{2}\rightarrow
\Lambda_{1}+\Lambda_{2}$. 
\qed
\end{proof}

Observe that if $\Lambda_2=\varnothing$, so we have a (uncoupled) block,
then $\BBB_{\Lambda_{1}+\Lambda_{2}}$ becomes the symmetric group $\SS_{\Lambda_1}$.
The following proposition summarizes all we need about the action of
$W(B_n)$ on $\HH(B_n)$:

\begin{proposition}\label{section:arrangementmonoids:result375}
In the action of the Weyl group $W(B_n)$ on $\HH(B_n)$, 
two subspaces $X(\Delta,\Lambda)$ and $X(\Delta',\Lambda')$ 
lie in the 
same orbit if and only if $|\Delta|=|\Delta'|$ and 
$\|\Lambda\|=\|\Lambda'\|$. The cardinality of the orbit of
the subspace $X(\Delta,\Lambda)$ is 
$$
2^{n-m-p-q}{n\choose n-m}\frac{(n-m)!}{b_\lambda},
$$
where $m=|\Delta|$, $\Lambda$ has the form
(\ref{section:arrangementmonoids:equation200}) and $b_\lambda$ is
given by (\ref{section:symmetricgroup:equation300}).
Moreover, if $W_X$ is the isotropy group
of $X=X(\Delta,\Lambda)$ then
\begin{equation}
  \label{section:symmetricgroup:equation500}
  W_X\cong\SS_{\pm\Delta}\times\BBB_{\Lambda_{11}+\Lambda_{12}}\times
\cdots
\times\BBB_{\Lambda_{q1}+\Lambda_{q2}}
\times\SS_{\Lambda_1}\times\cdots\times\SS_{\Lambda_p}
\end{equation}
\end{proposition}

The groups (\ref{section:symmetricgroup:equation500}) thus describe the
parabolics in $W(B_n)$, just as the Young subgroups do for $W(A_{n-1})$.

\begin{proof}
The orbit description and size is \cite{Orlik92}*{Proposition 6.75}. For the isotropy
group, we have $\xx g(\pi,T)=\xx$ for all $\xx\in X(\Delta,\Lambda)$ precisely when
$(\pi,T)\in\SS_n\ltimes\mathbf{2}^n$ leaves $\Delta$ and each block of
$\Lambda$ invariant. The expression 
for $W_X$ follows.
\qed
\end{proof}

In the following we identify the group
(\ref{section:symmetricgroup:equation500}) 
with a subgroup of $\SS_{\pm I}$. The proof uses 
Proposition \ref{section:arrangementmonoids:result375}, and is another application of
Theorem \ref{section:reflectionmonoids:result100} 
and (\ref{section:reflectionmonoids:equation200}):

\begin{theorem}\label{section:arrangementmonoids:result400}
The Coxeter arrangement monoid $M(B_n,\HH)$ has order
\begin{align*}
|M(B_n,\HH)|
=\sum_{(\Delta,\Lambda)}
[&\SS_{\pm I}:\SS_{\pm\Delta}\times\BBB_{\Lambda_{11}+\Lambda_{12}}\times
\cdots
\times\BBB_{\Lambda_{q1}+\Lambda_{q2}}
\times\SS_{\Lambda_1}\times\cdots\times\SS_{\Lambda_p}]\\
&=4^{n}(n!)^2\sum_{(m,\lambda)} \frac{1}{4^m(m!)^2d_\lambda},
\end{align*}
the first sum over all $(\Delta,\Lambda)\in\TT$, and the second over all 
pairs $(m,\lambda)$ where $0\leq m\leq n$ is an integer, 
$\lambda=(\lambda_1,\ldots,\lambda_p)$ is a partition
of $n-m$ and $d_\lambda=2^p\,b_\lambda\lambda_1!\ldots\lambda_p!$ with $b_\lambda$ as in
(\ref{section:symmetricgroup:equation300}).
\end{theorem}

\subsection{Coxeter arrangement monoids of type $D$}
\label{section:arrangementmonoids:typeD}

We repeat \S\ref{section:arrangementmonoids:typeB} for type $D$, just
running through the answers.
Replace $\TT$ by the sublattice $\TT^\circ$ consisting of the $(\Delta,\Lambda)$
with $|\Delta|\not=1$, the map $(\Delta,\Lambda)\mapsto X(\Delta,\Lambda)$
above restricting to an isomorphism $\TT^\circ\rightarrow\HH(D_n)$. 

Let $\mathbf{2}^I_+\subset\mathbf{2}^I$ be the subgroup consisting of those 
$T\subset I$ with $|T|$ even. Then the group isomorphism 
(\ref{section:arrangementmonoids:equation300})
restricts to an isomorphism $W(D_n)\rightarrow\SS_I\ltimes\mathbf{2}^I_+$.
The action (\ref{section:arrangementmonoids:equation400})
of $\SS_I\ltimes\mathbf{2}^I_+$ on $\TT^\circ$ models the action of
$W(D_n)$ on $\HH(D_n)$ as before.

\begin{proposition}
\label{section:arrangementmonoids:result500}
If $X(\Delta,\Lambda)$ and $X(\Delta',\Lambda')$ lie in the 
same orbit of the action of $W(D_n)$ on 
$\HH(D_n)$, then
$|\Delta|=|\Delta'|$ and $\|\Lambda\|=\|\Lambda'\|$. Conversely, suppose 
that $|\Delta|=|\Delta'|$ and $\|\Lambda\|=\|\Lambda'\|$.
\begin{enumerate}
\item If $|\Delta|\geq 2$ then $X(\Delta,\Lambda)$ and 
$X(\Delta',\Lambda')$ lie in the same orbit, which has cardinality 
as in Proposition \ref{section:arrangementmonoids:result375}.
\item If $\Delta=\varnothing$, then the $W(B_n)$ orbit 
determined by $\|\Lambda\|=(\lambda_{11}+\lambda_{12},
\ldots,\lambda_{q1}+\lambda_{q2},\lambda_1,\ldots,\lambda_p)$ 
forms a single $W(D_n)$ orbit, except when each $\lambda_{i1}+\lambda_{i2}$ 
and $\lambda_i$ are even, in which case it decomposes
into two $W(D_n)$ orbits of size
$$
\frac{2^{n-p-q-1}\,n!}{b_\lambda}.
$$
\item If $X=X(\Delta,\Lambda)\in\HH(D_n)\subset\HH(B_n)$, then
the isotropy groups $W(D_n)_X$ and $W(B_n)_X$
coincide when $\Delta=\varnothing$
and each $\lambda_{i1}+\lambda_{i2}$ and
$\lambda_i$ are even, otherwise, $W(D_n)_X$ has index $2$ in $W(B_n)_X$.
\end{enumerate}
\end{proposition}

\begin{proof}
The first two parts are just
\cite{Orlik92}*{Proposition 6.79}. For the third, 
the index of $W(D_n)_X$ in $W(B_n)_X$ is at most $2$ as $W(D_n)_X=W(D_n)\cap W(B_n)_X$ with 
$W(D_n)$ of index two in $W(B_n)$. Thus either $W(D_n)_X$ has index $2$ in $W(B_n)_X$
or the isotropy groups coincide, with the latter happening precisely when 
$Xg(\pi,T)=X$ for $g(\pi,T)\in W(B_n)$ implies that 
$g(\pi,T)$ is in $W(D_n)$, ie: that $|T|$ is even.
It is easy to check that this happens if and only if $\Delta=\varnothing$
and each $\lambda_{i1}+\lambda_{i2}$ and $\lambda_i$ is even.
\qed
\end{proof}

\begin{theorem}\label{section:arrangementmonoids:result600}
The Coxeter arrangement monoid $M(D_n,\HH)$ has order,
$$
|M(D_n,\HH)|=2^{2n-1}(n!)^2
\sum_{(m,\lambda)} \frac{\varepsilon_{m,\lambda}^{}}{4^{m}(m!)^2\,d_{\lambda}}
,
$$ 
the sum over all 
pairs $(m,\lambda)$ where $0\leq m\leq n$ is an integer $\not= 1$
and $\lambda=(\lambda_1,\ldots,\lambda_p)$ is a partition of $n-m$,
with $\varepsilon_{m,\lambda}^{}=1$ if $m=0$ and each $\lambda_i$ is even, and
$\varepsilon_{m,\lambda}^{}=2$ otherwise.
\end{theorem}

\subsection{Coxeter arrangement monoids of exceptional types}
\label{section:arrangementmonoids:exceptional}

Finally, to the orders of the Coxeter arrangement monoids in the exceptional cases,
where a combinatorial description
of $\HH$ is harder, but an enumeration of the orbits 
of the $W(\Phi)$-action on $\HH$--and their sizes and common isotropy
groups--suffices for our purposes. All this information is
contained in 
\cites{Orlik83,Orlik82} (see \cite{Orlik92}*{Appendix C}). For example, we 
can reproduce the essential information when $\Phi=F_4$ from
\cite{Orlik92}*{Table C.9, page 292} as
$$
A_0,12A_1,12\tilde{A}_1,72(A_1\times\tilde{A}_1),16A_2,16\tilde{A}_2,18B_2,
12C_3,12B_3,48(A_1\times\tilde{A}_2),48(\tilde{A}_1\times A_2),F_4,
$$
where each term $n\Phi$ indicates a $W(F_4)$-orbit on $\HH$
of size $n$ and with isotropy group $W(\Phi)$. For our purposes the tildes can be
ignored (so that $\tilde{A}_n=A_n$) and we also have $W(C_n)=W(B_n)$. The data can then be
plugged directly into (\ref{section:reflectionmonoids:equation200}),
using the orders given in Tables
\ref{section:reflectiongroups:table1}-\ref{section:reflectiongroups:table2} to get a calculation
for the order of the Coxeter arrangement monoid of type $F_4$,
$$
|M(F_4,\HH)|=
2^7\,3^2\biggl(1+\frac{12}{2}+\frac{12}{2}+\frac{72}{2^2}
+\cdots+\frac{48}{2^2\,3}+\frac{1}{2^7\,3^2}\biggr)=
11\cdot 4931.
$$

\begin{proposition}\label{section:arrangementmonoids:result700}
The orders of the exceptional Coxeter arrangement monoids are 
$$
\begin{tabular}{c|ccccc}
\hline
$\Phi$&$G_2$&$F_4$&$E_6$&$E_7$&$E_8$\\\hline
&&&&&\\
$|M(\Phi,\HH)|$
&$7^2$&$11\cdot 4931$&$2^4\cdot 5^2\cdot 40543$
&$3\cdot 113\cdot 24667553$
&$11\cdot 79\cdot 55099865069$\\
&&&&\\\hline
\end{tabular}
$$
\end{proposition}

What significance (if any) there is to these strange prime factorizations, we do not know.

\section{Generalities on reflection monoids}
\label{section:generalities}

We pause to explore some of the basic properties of reflection monoids.
Among other things, we justify the terminology ``reflection monoid''. 

Recall the definition of an inverse monoid $M$ from \S\ref{section:reflectionmonoids}. 
The sets of units and idempotents, 
$$
G=G(M)=\{a\in M\,|\, aa^{-1}=1=a^{-1}a\}\text{ and }
E=E(M)=\{a\in M\,|\, a^2=a\},
$$ 
form a subgroup 
and commutative submonoid respectively.
It is an elementary fact in semigroup theory that any 
commutative monoid of idempotents
carries the structure
of a meet semi-lattice, with order $e\leq f$ iff $ef=e$, and a unique
maximal element.
For this reason, $E$
is referred to as the \emph{semilattice of idempotents\/}.
An \emph{inverse submonoid\/} of an inverse monoid $M$ is simply a
subset $N$ that forms an inverse monoid under the same multiplicative
and $^{-1}$ operations. More details can be found in 
\cite{Howie95,Lawson98}.

We also observed in \S\ref{section:reflectionmonoids} that two primordial examples
of inverse monoids are the symmetric inverse monoid $\II_X$ and the 
general linear monoid $\ml(V)$, which have units the symmetric group $\SS_X$
and general linear group $GL(V)$ respectively. The idempotents $E(\II_X)$ consist
of the partial identities $\ve_Y^{}$ for $Y\subset X$
Similarly
the idempotents of $\ml(V)$ are the partial identities on subspaces of $V$.

We shall be particularly interested in factorizable inverse monoids:
monoids $M$ with $M = EG=GE$. 
Factorizability captures formally an idea
used informally in \S\S\ref{section:reflectionmonoids}-\ref{section:symmetricgroup}
and \S\ref{section:booleanmonoids}:
if $\aa\in M$ 
where $M$ is an inverse submonoid of $\II_X$, we have
$\aa\in EG$ if and only if $\aa$ is a restriction of a unit of
$M$. Similarly for $\ml(V)$.
In particular, $\II_n$ and $\ml(V)$, for $V$ finite dimensional, are factorizable,
but $\II_X$ for $X$ infinite is not (we remind the reader of our running 
assumption that $V$ is finite dimensional).

Of course this paper is about the monoids $M(G,\cS)$ of Definition 
\ref{section:reflectionmonoids:definition100}, where $g_X^{}$ is a unit precisely
when $X=V$ and $g_X^2=g_X^{}$ precisely when the 
restriction of $g$ to $X$ is the identity on $X$. Thus 
$G$ is the group of units and 
the idempotents $E$ are
the partial identities $\ve_X^{}$ for $X\in\cS$. If
$\cS$ is ordered by inclusion,
then $X\mapsto\ve_X^{}$ is an 
isomorphism of meet semi-lattices $\cS\rightarrow E$.

\begin{proposition}\label{section:generalities:result100}
Let $V$ be a finite dimensional vector space. Then $M\subset\ml(V)$ is a 
factorizable inverse submonoid if and only if $M=M(G,\cS)$ for
$G$ the group of units of $M$ and $\cS=\{\dom\,\aa\,|\,\aa\in M\}$.
\end{proposition}

\begin{proof}
We observed in \S\ref{section:reflectionmonoids} that
$M(G,\cS)$ is a fatorizable inverse submonoid of $\ml(V)$.
Conversely, if $M$ is a factorizable inverse submonoid, let $G$ be its 
group of units and let $\cS=\{\dom\,\aa\,|\,\aa\in M\}$. Then $V=\dom\,\ve$,
for $\ve:V\rightarrow V$ the identity map,
$\dom\,\aa\cap\dom\,\bb=\dom(\aa\aa^{-1}\bb\bb^{-1})$, and 
$(\dom\aa)g=\dom(g^{-1}\aa)$ for $g\in G$. Thus $\cS$ is a system of subspaces 
for $G$ in $V$, allowing us to form the monoid $M(G,\cS)$ of partial
isomorphisms. If $g\in G$ and $X=\dom\,\aa\in\cS$ then $g_X^{}=\aa\aa^{-1}g\in M$,
so that $M(G,\cS)$ is a factorizable inverse submonoid of $M$. 
Since $M(G,\cS)$ contains all the units and all the idempotents of $M$,
it follows that $M(G,\cS)=M$.
\qed
\end{proof}

A \emph{partial reflection\/} of $V$ is a partial isomorphism of the 
form $s_X^{}$ where $s$ is a full reflection.
If a reflection group is a group generated by reflections, then a reflection
monoid ought to be a monoid generated by partial reflections. In fact, a reflection
monoid is somewhat more--this is a consequence of the structure placed on the domains
through the system $\cS$:

\begin{corollary}\label{section:generalities:result200}
A submonoid $M\subset\ml(V)$ 
is a reflection monoid if and
only if $M$ is a factorizable inverse monoid generated by partial reflections.
\end{corollary}

\begin{proof}
If $M=M(W,\cS)$ for $W$ a reflection group, then $M$ is a factorizable inverse
submonoid of $\ml(V)$ with any
element having the form $\ve_X^{}g$ for some $X\in \cS$ and $g\in W$. Now 
$g =s_1\dots s_k$ for some reflections $s_1,\dots,s_k$ and $\ve_X^{}s_1$ is
a partial reflection, so $\ve_X^{}g = (\ve_X^{}s_1)s_2\dots s_k$ is
a product of partial reflections. Conversely, if $M$ is a 
factorizable inverse submonoid of $\ml(V)$ then $M=M(G,\cS)$ for $G=G(M)$
by Proposition \ref{section:generalities:result100}. 
It is easy to see that the non-units in $M$
form a subsemigroup, and hence every unit of $M$ must be a
product of (full) reflections, so $G$ is a
reflection group. Indeed, if $S$ is the set of generating partial
reflections for $M$, then $G=\langle S'\rangle$, for
$S'\subset S$ the full reflections.
\qed
\end{proof}

If $V_i$ $(i=1,2)$ are vector spaces with $G_i\subset GL(V_i)$ and 
$\cS_i$ are systems in $V_i$ for $G_i$, then a monoid homomorphism 
$M(G_1,\cS_1)\rightarrow M(G_2,\cS_2)$ induces a group 
homomorphism $G_1\rightarrow G_2$ and a homomorphism
$E_1\rightarrow E_2$ of semilattices. 
By the comments immediately prior to Proposition \ref{section:generalities:result100},
the latter is equivalent to
a poset map $\cS_1\rightarrow\cS_2$ between the two systems, ordered by inclusion.
If the homomorphism 
$M(G_1,\cS_1)\rightarrow M(G_2,\cS_2)$
is an isomorphism we have isomorphisms
of groups $G_1\rightarrow G_2$ and posets $\cS_1\rightarrow\cS_2$.

Recall that Green's relation $\mathscr{R}$ on a monoid $M$ is defined
by the rule that $a\mathscr{R} b$ if and only if $aM = bM$. The relation
$\mathscr{L}$ is the left-right dual of $\mathscr{R}$; we define 
$\mathscr{H} = \mathscr{R} \cap \mathscr{L}$ and $\mathscr{D} =
\mathscr{R}\ \vee \mathscr{L}$. In fact, 
$\mathscr{D} = \mathscr{R} \circ \mathscr{L} = \mathscr{L} \circ \mathscr{R}$.
Finally, $a\JJ b$ if and only if $MaM =
MbM$. In an inverse monoid, $a\mathscr{R} b$ if and only if $aa^{-1} =
bb^{-1}$ and similarly, $a\mathscr{L}b$ if and only if $a^{-1}a =
b^{-1}b$. 

\begin{proposition} \label{section:generalities:result300}
Let $\aa,\bb$ be elements of the monoid $M = M(G,\cS)$ of partial
linear isomorphisms, with
$\aa = g_X^{}$ and $\bb = h_Y^{}$ where $g,h \in G$ and $X,Y\in \cS$. Then
\begin{description}
\item[(i).] $\aa \mathscr{R} \bb$ if and only if $X = Y$;
\item[(ii).] $\aa \mathscr{L} \bb$ if and only if $Xg = Yh$;
\item[(iii).] $\aa \mathscr{D} \bb$ if and only if $Y\in XG$;
\item[(iv).] if $\cS$ consists of finite dimensional spaces, then
$\mathscr{J} = \mathscr{D}$.
\end{description}
\end{proposition}



Parts (i) and (ii) follow from \cite{Howie95}*{Proposition 2.4.2} and the well 
known characterization of $\mathscr{R}$ and $\mathscr{L}$ in $\ml(V)$. The rest is now 
a straightforward exercise for the reader.

\section{Linear algebraic monoids and Renner monoids}
\label{section:renner}

Given a linear algebraic group $\G$ one can extract from it
a finite group, the Weyl group, which turns out to play a number
of roles. On one hand, it acts as a group
of reflections of a space naturally associated to $\G$, and as such is a Weyl group in the sense
of \S\ref{section:reflectiongroups}. 
On another, there is the Bruhat decomposition of $\G$ with respect to a Borel subgroup, 
and the terms in the decomposition are parametrized by the elements of the Weyl group.

The theory of linear algebraic \emph{monoids\/} was developed independently, and then
subsequently collaboratively, by Mohan Putcha and Lex Renner, in the 1980's.
Among its chief achievements is the classification 
\cites{Renner85a,Renner85b} 
of the reductive monoids, and the formulation of a Bruhat decomposition
\cite{Renner86}, analogous to that for groups. 
%

In this section we show that one can also extract, in a natural fashion,
two finite monoids from a linear algebraic monoid $\M$. The first, which is new, is 
a reflection monoid in the same space that the Weyl group is a reflection group.
The other, which is well known and coined the \emph{Renner monoid\/}
by Solomon \cite{Solomon90}, plays the same role as the Weyl group in the 
Bruhat decomposition of $\M$. In general it does not seem to be possible to find a single monoid
to play all the roles that the Weyl group plays. Nevertheless our main result,
Theorem \ref{section:renner:result100} below, shows that
this reflection monoid and the Renner monoid are very closely related.

The prerequisites for this section are more demanding than for earlier ones, and
the reader who is unfamiliar with the theory of algebraic groups may find it helpful
at first to think in terms of Example \ref{section:renner:example100}. 
In any case, standard references on algebraic groups are 
\cites{Borel91,Humphreys75,Springer98}, and on algebraic monoids, the books of
Putcha and Renner \cites{Putcha88,Renner05}. We particularly recommend the
excellent survey of Solomon \cite{Solomon95}.

Throughout, $k$ is an algebraically closed field. 
An \emph{affine\/} (or \emph{linear\/}) \emph{algebraic monoid\/} 
$\M$ over $k$ is an affine algebraic variety 
together
with a morphism $\varphi:\M\times\M\rightarrow\M$ of varieties,
such that the product $xy:=\varphi(x,y)$ gives $\M$ the structure of a monoid
(ie: $\varphi$ is an associative morphism of varieties and there is
a two-sided unit $1\in\M$ for $\varphi$).
We will assume that the monoid $\M$ is \emph{irreducible\/}, that is, the
underlying variety is irreducible, in which case the group $\G$ of units 
is a connected algebraic group with $\overline{\G}=\M$ (Zariski closure).
Adjectives normally applied to $\G$ are then transferred to $\M$; thus we have
reductive monoids, simply connected monoids, soluble monoids, and so on.
We have, in analogy to the group case, that any affine algebraic monoid can be 
embedded as a closed submonoid in $M_n(k)$ for some $n$.

From now on, let $\M$ be reductive.
The key players, just as they are for algebraic groups, 
are the maximal tori $T\subset\G$ and their closures
$\overline{T}\subset\M$. Let $\XXX(T)$ be the character group of all
morphisms of algebraic groups $\chi:T\rightarrow\G_m$ (with 
$\G_m$ the multiplicative group of $k$) and $\XXX(\overline{T})$
similarly the commutative monoid of morphisms of $\ov{T}$ into $k$ as
a multiplicative monoid. Then $\XXX(T)$
is a free $\Z$-module, and restriction (together with the denseness
of $T$ in $\ov{T}$) embeds $\XXX(\ov{T})\hookrightarrow\XXX(T)$.

The Weyl group $W_\G=N_\G(T)/T$ of automorphisms of $T$ acts faithfully on 
$\XXX(T)$ via $\chi^g(t)=\chi(g^{-1}tg)$, thus realizing an injection
$W_\G\hookrightarrow GL(V)$ for $V=\XXX(T)\otimes\R$.
We will abuse notation and write $W_\G$ for both the Weyl group and its image in $GL(V)$. 
The non-zero weights $\Phi:=\Phi(\G,T)$ of the adjoint representation
$\G\rightarrow GL(\goth{g})$ form a root system in $V$ with 
the Weyl group $W_\G$ the reflection group $W(\Phi)$ associated to $\Phi$
(here, $\goth{g}$ is the associated Lie algebra).

We now need a digression
to review some basic facts about convex polyhedral cones, 
for which we follow \cite{Fulton93}*{\S 1.2}. If $V$ is a real space and $\vv_1,\ldots,\vv_s$
a finite set of vectors, then the convex polyhedral cone with generators $\{\vv_i\}$
is the set $\ss=\sum\lambda_i\vv_i$ with $\lambda_i\geq 0$. The dual cone 
$\ss^\vee\subset V^*$ consists of those $\uu\in V^*$ taking non-negative values
on $\ss$. A face $\tau\subset\ss$ is the intersection with $\ss$ of the kernel $\uu^\perp$
of a $\uu\in\ss^\vee$, and the faces form a meet semilattice $\FF(\ss)$ under inclusion.
If $\tau\in\FF(\ss)$ 
and $\ov{\tau}$ is the $\R$-span in $V$ of $\tau$, 
then $\ss\cap\ov{\tau}=\tau$.
In particular, if 
$\bigcap\ov{\tau}_i=\bigcap\ov{\mu}_j$ in $V$ then 
$\bigcap\tau_i=\bigcap\mu_j$ in $\FF(\ss)$.

\begin{figure}
  \centering
\begin{pspicture}(0,0)(14.3,4.5)
\rput(7.15,2.2){\BoxedEPSF{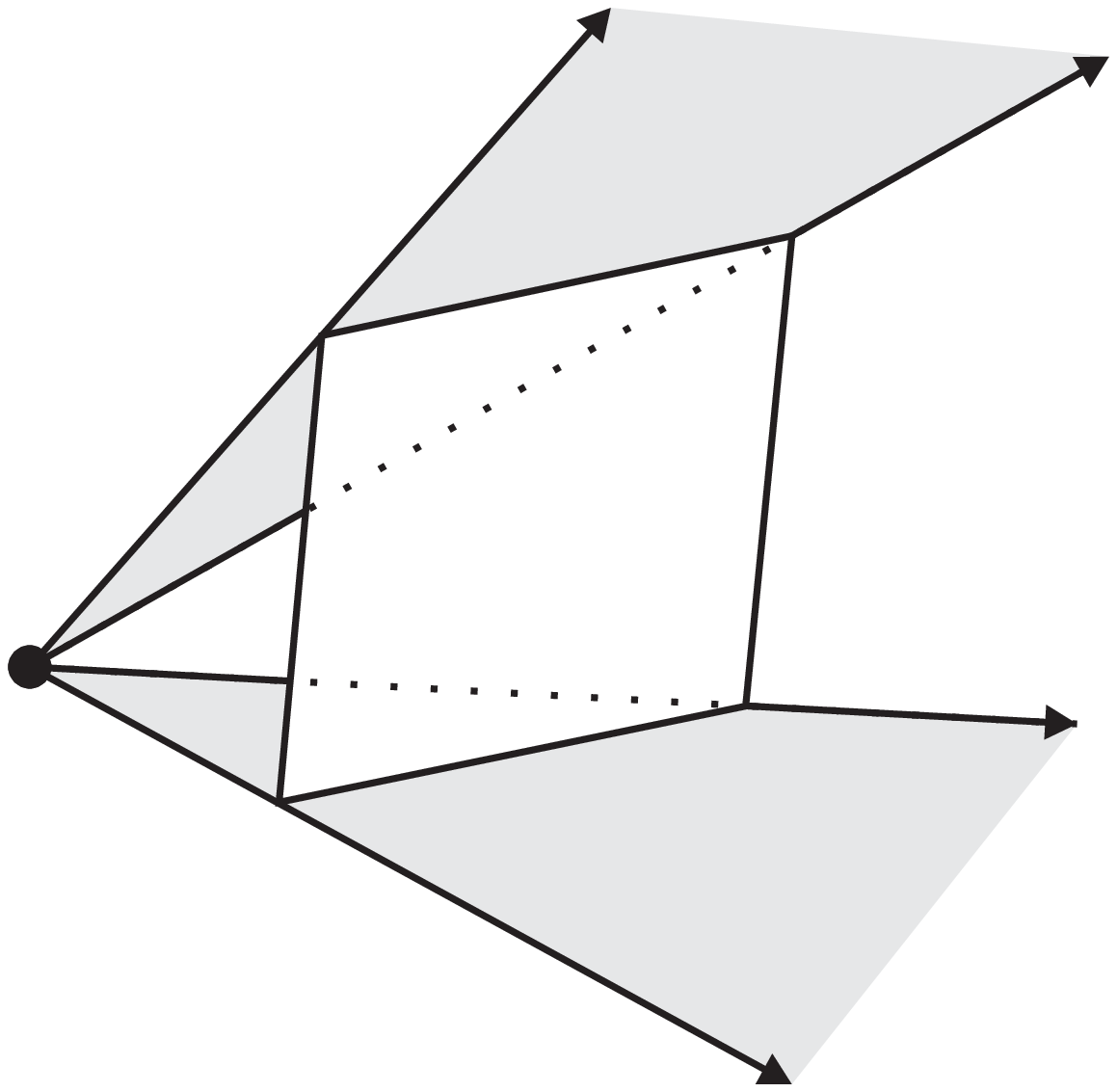 scaled 400}}
\rput(8.5,2.5){$\ss$}\rput(8,4){$\tau_1$}\rput(7.5,1){$\tau_2$}
\rput(5.95,3.2){$\vv_1$}\rput(8.4,3.4){$\vv_2$}
\rput(8.2,1.7){$\vv_3$}\rput(6,.85){$\vv_4$}
\end{pspicture}
\caption{Failure of $\ov{\tau}\subset\bigcap\ov{\tau}_j$ to be an equality.}
\label{section:renner:fig100} 
\end{figure}

If $\{\tau_j\}\subset\FF(\ss)$ are faces of $\ss$ and $\tau=\bigcap\tau_j$, then we have
$\ov{\tau}\subset\bigcap\ov{\tau}_j$. In general this inclusion is not an equality,
a simple observation with non-trivial consequences. The reader interested
in the source of the failure of the homomorphism of Theorem \ref{section:renner:result100} 
to be an
isomorphism can trace it back to here. The cone in Figure \ref{section:renner:fig100} has
faces $\tau_1,\tau_2$ with $\ov{\tau}$ 
the zero subspace but $\bigcap\ov{\tau}_j$ $1$-dimensional.

A cone is simplicial if it has a set $A=\{\vv_i\}$ of linearly 
independent generators. 
If $\tau_i$ is the cone on $\{\vv_1,\ldots,\wh{\vv}_i,\ldots,\vv_s\}$,
then $\tau_i=\ss\cap\uu_i^\perp$, where $\uu_i$ is the vector corresponding to $\vv_i$ in the
dual basis for $V^*$. Thus $\tau_i$ is a face of $\ss$, and the face lattice
$\FF(\ss)$ is isomorphic to the lattice of all subsets of $A$--or, if one prefers, to 
the Boolean lattice on the $1$-dimensional
faces $\R^+\vv_i$ of $\ss$.
If $\tau\in\FF(\ss)$ corresponds to $A_\tau\subset A$ then 
$\tau_1\cap\tau_2$ corresponds to $A_{\tau_1}\cap A_{\tau_2}$, and
$\ov{\tau}$ is the $\R$-span of  $A_\tau$.
In particular,
the $\R$-span of $\bigcap A_{\tau_j}$ is the intersection of the 
$\R$-spans of the  $A_{\tau_j}$,
and so if $\tau=\bigcap\tau_j$, we have $\ov{\tau}=\bigcap\ov{\tau}_j$ 
when $\ss$ is simplicial.
Finally, a cone is strongly convex if the dual $\ss^\vee$ spans $V^*$. Simplicial
cones are strongly convex. On the other hand, if $\dim V=2$, then
any strongly convex cone is simplicial \cite{Fulton93}*{1.2.13}.

Returning to our algebraic monoid
we may assume, by conjugating suitably, that the maximal
torus $T$ is a subgroup of the group
$\T_n$ of invertible diagonal matrices, where $n$ is the rank of $\G$. 

\begin{definition}\label{section:renner:definition100} 
Let $\M$ be a reductive algebraic monoid.
If $\chi_j$ is the restriction to $T$ of the
$j$-th coordinate function on $\T_n$, then let $\ss\subset\XXX(T)\otimes\R$ be the cone
given by $\ss=\sum\R^+\chi_i$.
\end{definition}

This cone will turn out to have a number of nice properties,
the first of which being that 
the character monoid $\XXX(\overline{T})=\ss\cap\XXX(T)$.
Secondly, 
the Weyl group $W_\G$, in its reflectional action on
$V$, acts on $\ss$, and this induces an action 
$\tau\mapsto\tau g$ of $W_\G$ on the face lattice $\FF(\ss)$.

Finally, and most importantly, 
the face lattice $\FF(\ss)$ models idempotents:
there is a lattice 
isomorphism 
$\FF(\ss)\rightarrow E(\ov{T})$, written $\tau\mapsto e_\tau$,
that is  $W_\G$-equivariant
with respect to the Weyl group actions on $\FF(\ss)$ and $E(\ov{T})$.
In short, $e_\tau^g=e_{\tau g}$ for any $\tau\in\FF(\ss)$
and $g\in W_\G$.
Solomon \cite{Solomon95}*{Corollary 5.5}, working instead with the dual
cone $\ss^\vee$ in the group of $1$-parameter subgroups of $T$,
has a lattice \emph{anti}-isomorphism $\FF(\ss^\vee)\rightarrow E(\ov{T})$.

\begin{example}\label{section:renner:example100} 
Let $\M=M_n(k)$, the monoid of $n\times n$ matrices over $k$, which can be naturally 
identified with $n^2$-dimensional affine space over $k$. The units are the
general linear group $\G=GL_n(k)$, a reductive group. The tori in $\G$ are the conjugates of
subgroups of diagonal matrices, and an
example of a maximal torus $T$ is the subgroup of all
invertible diagonal matrices,
$$
\begin{pspicture}(0,0)(14.3,2)
\rput(7.15,1){
$\text{diag}(\lambda_1,\ldots,\lambda_n):=
\left[
  \begin{array}{ccc}
  \lambda_1&&\\
  &\vrule width 0 mm height 0 mm depth 4 mm
   \vrule width 4 mm height 0 mm depth 0 mm&\\
  &&\lambda_n\\    
  \end{array}
\right]
\text{ with }\lambda_1\ldots\lambda_n\not=0.$
}
\rput(6.85,.65){\psline[linewidth=.4mm,linestyle=dotted]{-}(0,.6)(.6,0)}
\rput(6.8,.6){$\text{\Large{0}}$}\rput(7.5,1.5){$\text{\Large{0}}$}
\end{pspicture}
$$
Any other maximal torus is a conjugate of this one. The Zariski closure $\ov{T}$ consists
of all the diagonal matrices with no restriction on the $\lambda_i$. The map
$\chi_i:T\rightarrow \G_m$ sending $\text{diag}(\lambda_1,\ldots,\lambda_n)$ to 
$\lambda_i$ is a character of $T$, with an arbitrary character $\chi\in\XXX(T)$
having the form $\chi=\chi_1^{t_1}\ldots\chi_n^{t_n}$ for some $t_i\in\Z$. Thus
$\XXX(T)$ is a free $\Z$-module with basis the $\chi_i$ and $V=\XXX(T)\otimes\R$
is an $n$-dimensional space.

The normalizer $N_\G(T)$ consists of the monomial matrices: those 
having precisely one non-zero entry in each row and column. The Weyl group
$W_\G=N_\G(T)/T$ thus consists of the permutation matrices 
$A(\pi)=\sum_{i=1}^nE_{i,i\pi}$, where $E_{ij}$ is the matrix with
$(i,j)$-th entry $1$ and all other entries $0$, and $\pi\in\SS_n$.
The $W_\G$-action $\chi^g(t)=\chi(g^{-1}tg)$ 
on $V$ becomes the $\chi_i^{A(\pi)}=\chi_{i\pi}$ of \S\ref{section:symmetricgroup}, with 
the weights $\Phi=\{\chi_i\chi_j^{-1}\,|\,1\leq i\not= j\leq n\}$ the root
system $A_{n-1}$ of \S\ref{section:reflectiongroups}.

The idempotents $E(\ov{T})$ are those $\text{diag}(\lambda_1,\ldots,\lambda_n)$ 
with $\lambda_i\in\{0,1\}$ for all $i$. Indeed, the idempotents are parametrized
by the subsets of $\{1,\ldots,n\}$ according to the positions on the diagonal
of the $1$'s. Moreover, the order $e\leq f\Leftrightarrow ef=e$ corresponds
precisely to the subset for $e$ being contained in the subset for $f$. Thus
$E(\ov{T})$ is isomorphic as a lattice to the Boolean lattice of all subsets
of $\{1,\ldots,n\}$ with unique minimal element the zero matrix and
unique maximal element the identity
matrix.

The cone $\ss=\sum\R^+\chi_i$ is the positive quadrant in $V$, hence simplicial,
and with the integral points in $\ss$
the characters in the monoid $\XXX(\ov{T})$. The $W_\G$-action on $\ss$
permutes the vertices of the $(n-1)$-simplex consisting of the points
$\sum\mu_i\chi_i$ with $\sum\mu_i=1$, and the face lattice $\FF(\ss)$ can be
identified with the face lattice of this simplex.
Figure \ref{section:renner:figure200}
illustrates it all when $n=3$ with the $2$-simplex shaded and 
$\text{diag}(\lambda_1,\ldots,\lambda_n)$ represented by
$\lambda_1\ldots\lambda_n$.

\begin{figure}[h]
  \centering
\begin{pspicture}(0,0)(14.3,4)
\rput(0,-.2){
\rput(7,2.4){\BoxedEPSF{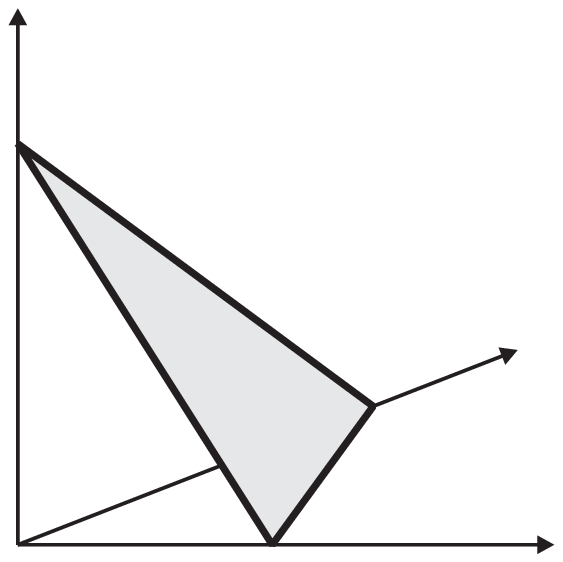 scaled 700}}
\rput(7,.3){$100$}\rput(8,1.4){$010$}\rput(4.7,3.4){$001$}
\rput(5.8,1.6){$101$}\rput(7,2.5){$011$}\rput(7.7,.9){$110$}
}
\end{pspicture}
  \caption{The simplicial cone $\ss\subset\XXX(T)\otimes\R$ and the 
isomorphism $\FF(\ss)\rightarrow E(\ov{T})$ for $\M=M_n k$.}
  \label{section:renner:figure200} 
\end{figure}

Finally, an abstract simplex with vertices $X=\{\chi_1,\ldots,\chi_n\}$ has faces 
the subsets of $X$, with inclusion of faces corresponding to inclusion of subsets. Thus the 
face lattice $\FF(\ss)$ is the Boolean lattice of subsets of $X$ and 
the isomorphism $\FF(\ss)\rightarrow E(\ov{T})$ sends $\{\chi_{i_1},\ldots,\chi_{i_k}\}$
to the diagonal matrix with $1$'s in positions
$\{i_1,\ldots,i_k\}$ and $0$'s elsewhere.
\end{example}

Returning to generalities, the point of the cone $\ss$ 
is that it gives a system of subspaces in $V$ for the
Weyl group:

\begin{definition}\label{section:renner:definition200}
Let $\M$ be a reductive algebraic monoid with $\ss$ the cone of
Definition \ref{section:renner:definition100} and $\FF(\ss)$ its face lattice.
Define the \emph{system $\cS_\M$ for $W_\G$ in $V=\XXX(T)\otimes\R$\/} 
to be the set of all intersections 
$\bigcap\ov{\tau}_j$, for $\tau_j\in\FF(\ss)$, and with the empty intersection taken to
be $V$.
\end{definition}

The terminology presupposes that $\cS_\M$ is indeed a system for $W_\G$,
but this is immediate as $W_\G$ acts on the face lattice $\FF(\ss)$.
One can also describe $\cS_\M$ as 
the system for $W_\G$ \emph{generated\/} by the subspaces $\ov{\tau}$.

Revisiting Example \ref{section:renner:example100}, the system $\cS_\M$ is just the
Boolean system of \S\ref{section:symmetricgroup} consisting of all $X(J)=\sum_{j\in J}\R\chi_j$
for $J\subset\{1,\ldots,n\}$.

\begin{definition}\label{section:renner:definition300} 
Let $\M$ be a reductive algebraic monoid with group of units $\G$.
The \emph{reflection monoid associated to $\M$\/}
is the monoid $M(W_\G,\cS_\M)$ where $W_\G$ is the Weyl group of $\G$ and
$\cS_\M$ is the system for $W_\G$ given in Definition \ref{section:renner:definition200}.
\end{definition}

Wrapping up Example \ref{section:renner:example100}, the Weyl group $W_\G$ acts
on $V$ as the symetric group $\SS_n$ permuting coordinates, and the system $\cS_\M$
is the Boolean one: the reflection monoid associated to $\M$ is thus our first example
$M(\SS_n,\BB)$ from \S\ref{section:symmetricgroup}--the
symmetric inverse monoid $\II_n$ or the Boolean reflection monoid $M(A_{n-1},\BB)$.

We now introduce a second intermediate monoid associated with $\M$
that is constructed via
the cone $\ss$ and uses the idea of a system of \emph{subsets\/}
for $W_\G$; this idea is defined in general in
\S\ref{section:complements:factorizable}. For now, we observe that
$\ss\in\FF(\ss)$, the lattice $\FF(\ss)$ is closed under intersection
and $W_\G$ acts on $\FF(\ss)$. These facts allow us
to form a inverse monoid $M(W_\G,\FF(\ss))$ in the same way that
$M(G,\cS)$ is constructed when $\cS$ is a system of subspaces for a
group $G$. The elements of $M(W_\G,\FF(\ss))$ have the form
$w_\tau, (w\in W_\G,\tau\in\FF(\ss))$, where $w_\tau$ is the restriction
of $w$ to $\tau$. Multiplication is given by
$w_{\tau}w'_{\tau'}=(ww')_{\mu}$ where $\mu=\tau\cap\tau'w^{-1}$.

The minimum face of $\ss$ is $\ss\cap-\ss$, which is the largest subspace of $V$
contained in $\ss$. Write $Z$ for this subspace and let
$\pi:V\rightarrow V/Z$ be the canonical homomorphism. Then
$\pi(\ss)$ is a strongly convex polyhedral cone in $V/Z$,
and the face lattices $\FF(\ss)$ and $\FF(\pi(\ss))$ are isomorphic.

Not suprisingly, there is a close connection between the reflection
monoid $M(W_\G,\cS_\M)$
and $M(W_\G,\FF(\ss))$. First, note that
idempotents in $M(W_\G,\cS)$ are products $\ve_X^{}=\prod\ve_j$
where $X=\bigcap\ov{\tau}_j\in\cS_\M$ and $\ve_j$ is the partial
identity on $\ov{\tau}_j$. Factorizability (Proposition
\ref{section:generalities:result100}) thus gives that any element of the
reflection monoid has the form $\ve_X^{} w=(\prod\ve_j)w$ for $w\in W_\G$. 

Define $\theta:M(W_\G,\cS_\M)\rightarrow M(W_\G,\FF(\ss))$ by 
$$
(\ve_X^{} w)\theta=\left(\prod e_j\right)w,
$$
where $e_j$ is the partial identity $\ve_{\tau_j}$ on
$\tau_j$.

\begin{proposition}\label{section:renner:result100} 
The map $\theta:M(W_\G,\cS_\M)\rightarrow M(W_\G,\FF(\ss))$ is a
surjective homomorphism, and is an isomorphism if 
and only if the cone $\pi(\ss)$ is simplicial.
\end{proposition}

\begin{proof}
The identity map of $W_\G$ is an isomorphism between the groups of
units of the two monoids. Let $\varphi$ be the map between the respective
idempotents given by $(\prod\ve_j)\varphi=\prod e_j$. 
As 
$\bigcap\ov{\tau}_i=\bigcap\ov{\mu}_j$ implies
$\bigcap\tau_i=\bigcap\mu_j$ in $\FF(\ss)$, the map $\varphi$ is well
defined. It is clearly a homomorphism and surjective.
Next, note that $\ov{\tau}_jw=\ov{\tau_jw}\,(w\in W_\G)$, and so
$$
w\ve_X^{}w^{-1}
=\prod\ve_{\ov{\nu}_j}^{},\,\,(\ov{\nu}_j=\ov{\tau_j w}),
$$
giving
$$
(w\ve_X^{}w^{-1})\varphi=
w^{-1}\left(\prod\ve_{\tau_j}\right)w
=w^{-1}\left(\prod\ve_j\right)\kern-1mm\varphi\,w.
$$
Finally, $(\prod\ve_j)w=\prod\ve_j$
if and only if $w$ leaves $\bigcap\ov{\tau}_j$ fixed pointwise, and
this implies that $w$ leaves $\bigcap\tau_j$ fixed pointwise, so
$(\prod e_j)w=\prod e_j$.

We thus have all the ingredients needed to apply Proposition
\ref{section:reflectionmonoids:result200} and get $\theta$ a
surjective homomorphism. 

Now suppose that $\pi(\ss)$ is not simplicial. Then $\dim\pi(\ss)>2$
so $\dim\ss>2+\dim Z$. There are maximal faces in $\FF(\pi(\ss))$
which intersect in $\{Z\}$ and corresponding to maximal faces
$\tau_1,\tau_2$ in $\FF(\ss)$ with $\tau_1\cap\tau_2=Z$. As the spaces
$\ov{\tau}_i\,(i=1,2)$ are hyperplanes in $\ov{\ss}$, the intersection
$\ov{\tau}_1\cap\ov{\tau}_2$ has codimension two in $\ov{\ss}$ and hence $Z$ is
strictly contained in $\ov{\tau}_1\cap\ov{\tau}_2$. This translates
into $\ve_1\ve_2\not=0$ in $M(W_\G,\cS_\M)$ but $e_1e_2=0$ in
$M(W_\G,\FF(\ss))$, where, as usual, $\ve_i$ is the partial identity
on $\ov{\tau}_i$ and $e_i$ is the partial identity on $\tau_i$. Thus,
$\theta$ fails to be injective even on the idempotents of $M(W_\G,\cS_\M)$,
and so on $M(W_\G,\cS_\M)$ itself. 

On the other hand suppose that $\pi(\ss)$ is simplicial. Let
$X,Y\in\cS_\M$ with $\ve_X^{}\varphi=\ve_Y^{}\varphi$. If
$X=\bigcap\ov{\tau}_j$ and $Y=\bigcap\ov{\mu}_i$, then
$\ve_X^{}=\prod\ve_j$ and $\ve_Y^{}=\prod\ve'_i$, $\ve'_i$ the partial
identity on $\ov{\mu}_i$. Thus $\prod e_j=\prod e'_i$ where $e'_i$ is
the partial identity on $\mu_i$. This is equivalent to
$\bigcap\tau_j=\bigcap\mu_i$ and so in $\pi(\ss)$ we have
$\bigcap\pi(\tau_j)=\bigcap\pi(\mu_i)$. Since $\pi(\ss)$ is
simplicial, this gives that 
$$
\bigcap\ov{\pi(\tau_j)}=\ov{\bigcap\pi(\tau_j)}
=\ov{\bigcap\pi(\mu_i)}=\bigcap\ov{\pi(\mu_i)}.
$$ 
Now $Z\subset\tau_j$ for all $j$, so
$\bigcap\ov{\pi(\tau_j)}=\pi\left(\bigcap\ov{\tau}_j\right)$ and similarly 
$\bigcap\ov{\pi(\mu_i)}=\pi\left(\bigcap\ov{\mu}_i\right)$, whence
$\bigcap\ov{\tau}_j=\bigcap\ov{\mu}_i$, and so 
$\ve_X^{}=\ve_Y^{}$. Thus $\varphi$ is an isomorphism.

To apply Proposition \ref{section:reflectionmonoids:result200}, we
also need to show that $(\prod e_j)w=\prod e_j$ implies
$(\prod\ve_j)w=\prod\ve_j$. The assumption is equivalent to $w$
leaving $\bigcap\tau_j$ fixed pointwise, and hence $w$ leaves the
$\R$-span of $\bigcap\tau_j$ fixed pointwise. But this subspace is
equal to $\bigcap\ov{\tau}_j$ as $\ss$ is simplicial, and so
$(\prod\ve_j)w=\prod\ve_j$ follows. Thus $\theta$ is an
isomorphism by Proposition
\ref{section:reflectionmonoids:result200}.
\qed
\end{proof}

Now to monoid number three, 
where we can be briefer. The \emph{Renner monoid\/}
$R_\M$ of $\M$ is defined to be
$R_\M=\overline{N_\G(T)}/T$.
See \cite{Renner86}
or \cite{Putcha88}*{Chapter 11}.

Just as $\II_n$ is the archetypal inverse monoid, and as 
$M(A_{n-1},\BB)$ with $\BB$ the Boolean system it is the archetypal reflection monoid,
so in its incarnation as the \emph{rook monoid\/} it is the standard example of a Renner
monoid, namely for $\M=M_n(k)$ in Example 
\ref{section:renner:example100} 
above. The elements of the rook monoid are
the $n\times n$ matrices having $0,1$ entries with at most one non-zero entry
in each row and column. The etymology of ``rook monoid'' is that
each element represents an $n\times n$
chessboard with the $0$ squares empty, the $1$ squares containing
rooks and
the rooks mutually non-attacking.
The Renner monoids have also been explicitly described
in some other cases, for example when $\M$ 
is the ``symplectic monoid'' $MSp_n(k)=\ov{k^*Sp_n(k)}\subset M_n(k)$
\cite{Renner03}.

In general, we have $E(R_\M)=E(\ov{T})$ and $R_\M$ is a (factorizable)
inverse monoid. Consider, for $E=E(R_\M)$ the fundamental representation
$\aa:R_\M\rightarrow\II_E$, written $a\aa=\aa_a$, where
$\aa_a:Eaa^{-1}\rightarrow Ea^{-1}a$ is defined by $x\aa_a=a^{-1}xa$
for all $x\in Eaa^{-1}$. 
We shall describe the fundamental representation in a more general
context, in particular being more precise about the location of its
image, in \S\ref{section:complements:fundamental}. For now, we record
that $\aa$ is a homomorphism \cite{Howie95}*{Theorem 5.4.4}.

In general, $\aa$ is not an isomorphism (see \cite{Renner05}*{Proposition 8.3}), 
but if the reductive monoid $\M$ has a $0$, then
it follows
from \cite{Putcha88}*{Proposition 11.1} and \cite{Howie95}*{Theorem
  5.4.4} that $\aa$ is an isomorphism. Further, by
\cite{Putcha88}*{Theorem 6.20}, the length of a maximal chain in
$E(\ov{T})$ is equal to $\dim T$, where $\dim T=m$ when
$T=\T_m$. From \cite{Solomon95}, $\dim T=\text{rank}\,\XXX(T)=\dim V$
for $V=\XXX(T)\otimes\R$. Recall that there is an isomorphism
$\FF(\ss)\rightarrow E(\ov{T})$ which sends $\tau$ to an idempotent
$e_\tau$. We conclude that $\FF(\ss)$ has a chain of length $\dim
V$, and since $\ov{\ss}=V$, it follows that the least member of such a
chain must be the trivial subspace $\mathbf{0}$. Hence $\ss$ is
strongly convex. 

Restricting $\aa$ to $W_\G$ gives an isomorphism $\aa:W_\G\rightarrow
W_\G\aa$. Let $E=E(\ov{T})$ and
$E'$ be the partial identities $\{\ve_\tau\,:\,\tau\in\FF(\ss)\}$. 
Then we have an isomorphism
$\bb:E'\rightarrow E\aa$ given by 
$$
\ve_\tau\bb=e_\tau\aa,
$$ 
where in turn $e_\tau \aa=e_{\tau\aa}$ is the identity on $Ee_\tau$. 
Now define $\varphi:M(W_\G,\FF(\ss))\rightarrow R_\M$ by $(\ve_\tau
w)\varphi=(\ve_\tau\bb)(w\aa)$.

\begin{proposition}\label{section:renner:result200} 
Let $\M$ be a reductive monoid with $0$. Then the map 
$\varphi:M(W_\G,\FF(\ss))\rightarrow R_\M$ is an isomorphism.
\end{proposition}

\begin{proof}
We show first that $\bb$ is equivariant. For $w\in W_\G$, $\tau\in\FF(\ss)$: 
\begin{align*}
  \label{eq:7}
(w^{-1}\ve_\tau w)\bb
=(w^{-1}\cdot w|_\tau)\bb
&=((w^{-1}w)|_{\tau w})\bb 
=(\ve_{\tau w})\bb
=e_{\tau w}\aa\\
&=(w^{-1}e_\tau w)\aa
=(w\aa)^{-1}(e_\tau \aa)(w\aa)
=(w\aa)^{-1}(\ve_\tau\bb)(w\aa),
\end{align*}
as required. If $\ve_\tau w=\ve_\tau$, then $w$ leaves $\tau$ fixed
pointwise, so that $\kappa w=\kappa$ for all $\kappa\leq\tau)$. Now 
$(\ve_\tau\bb)(w \aa)=\id_{Ee_\tau}\aa_w=\aa_w|_{Ee_\tau}$. Also,
$e\in Ee_\tau$ if and only if $e=e_\kappa$ for some
$\kappa\leq\tau$. Thus, for all $\kappa\leq\tau$,
$$
e_\kappa\aa_w=w^{-1}e_\kappa w=e_{\kappa w}=e_\kappa,
$$
so $\aa_w|_{Ee_\tau}$ is the identity on $Ee_\tau$ as required. Thus Proposition
\ref{section:reflectionmonoids:result200} gives $\varphi$ is a
surjective homomorphism. To see that $\varphi$ is an isomorphism, all
that remains is to show that $\aa_w|_{Ee_\tau}=\id_{Ee_\tau}$ implies
$\ve_\tau w=\ve_\tau$, that is, $w$ leaves $\tau$ fixed
pointwise. Since $\ss$ is strongly convex, $\tau$ contains one
dimensional faces. Let $\kappa$ be one such. Then $\kappa$ is a ray
with $e_\kappa\leq e_\tau$, so $e_\kappa\in Ee_\tau$ and hence
$e_\kappa\aa_w=e_\kappa$, that is, $e_{\kappa w}=e_\kappa$. Thus
$\kappa w=\kappa$. As $w$ acts on $V$ as a reflection, and so is
orthogonal, it leaves $\kappa$ fixed pointwise. Since this is so 
for all one dimensional faces contained in $\tau$, it follows that $w$
leaves $\tau$ fixed pointwise and $\varphi$ is an isomorphism.
\qed
\end{proof}

The main result of the section now follows from the preceding two Propositions:

\begin{theorem}\label{section:renner:result100} 
Let $\M$ be a reductive algebraic monoid with $0$. 
\begin{itemize}
\item Let $\G$ be its group
of units with $T\subset\G$ a maximal torus, $\XXX(T)$ the character group and
$W_\G$ the Weyl group;
\item Let $\ss\subset \XXX(T)\otimes\R$ be the polyhedral cone of 
Definition \ref{section:renner:definition100}, $\FF(\ss)$ its face
lattice, and 
$\cS_\M$ the system for $W_\G$
of Definition \ref{section:renner:definition200};
\item Let 
$M(W_\G,\cS_\M)$ be the reflection monoid associated to $\M$ and
$M(W_\G,\FF(\ss))$ the monoid given by the system of subsets $\FF(\ss)$;
\item Finally, let $R_\M$ be the Renner monoid of $\M$.
\end{itemize}
Then $R_\M\cong M(W_\G,\FF(\ss))$ and there is a surjective
homomorphism $M(W_\G,\cS_\M)\rightarrow R_\M$ 
which is an isomorphism if and only if $\ss$ is a simplicial cone.
\end{theorem}

\begin{example}\label{section:renner:example200} 
As an illustration of the lack of injectivity of $f$, 
let $\M$ be the (normalization of) the closure of $\text{Ad}(\G)k^*$
for $\G$ the adjoint simple group of type $B_2$. Then
\cite{Renner85a}*{Example 3.8.3}, $\dim(\XXX(T)\otimes\R)=3$
with $\ss$ a cone on a square 
(see \cite{Renner85a}*{Figure 6} or Figure \ref{section:renner:fig100}). 
If $\tau_i$, $(i=1,2)$ are 
the cones on opposite, non-intersecting faces of the square, then
$\tau_1\cap\tau_2=\{0\}$, whereas $\ov{\tau}_1\cap\ov{\tau}_2$ is 
a $1$-dimensional subspace. 

\begin{figure}[h]
  \centering
\begin{pspicture}(0,0)(15,4)
\rput(-.5,0){
\rput(5.25,2){\BoxedEPSF{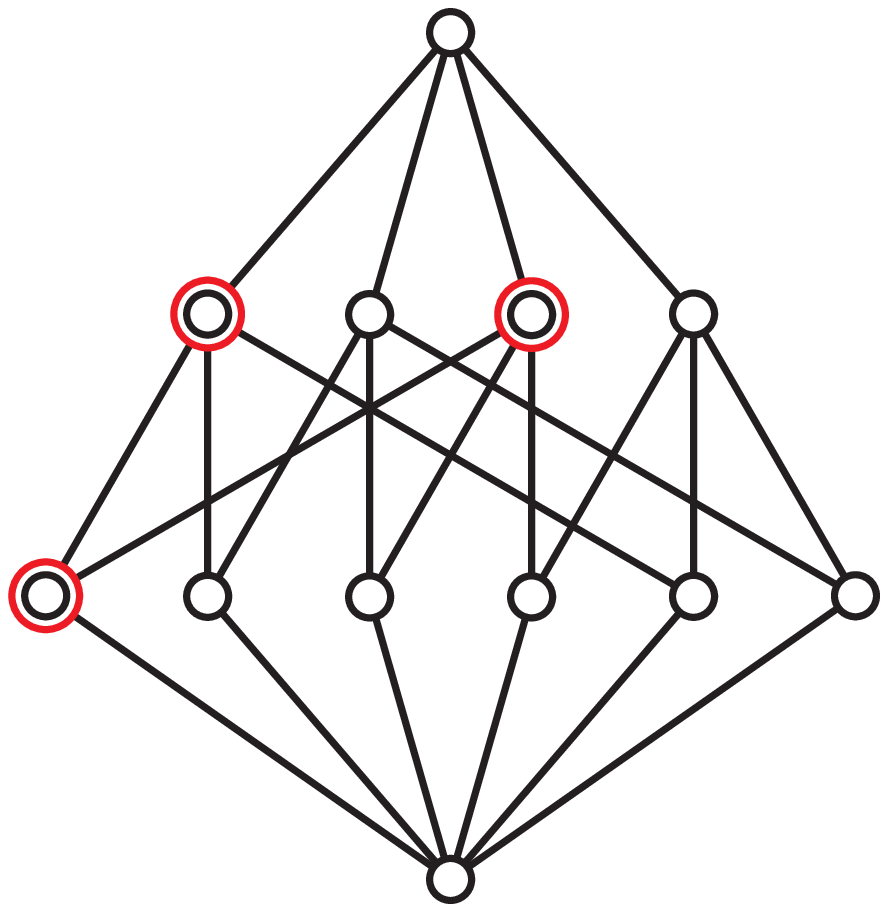 scaled 400}}
\rput(4,3.4){$\cS$}
\rput(3.9,2.6){{\red $\ov{\tau}_1$}}\rput(5.95,2.6){{\red $\ov{\tau}_2$}}
\rput(2.9,1.4){{\red $\ov{\tau}_1\cap\ov{\tau}_2$}}
}
\rput(7.35,2.25){$f$}
\psline[linewidth=.3mm]{->}(6.75,2)(7.95,2)
\rput(0,0){
\rput(10.6,3.4){$E(R_\M)$}
\rput(9.25,2){\BoxedEPSF{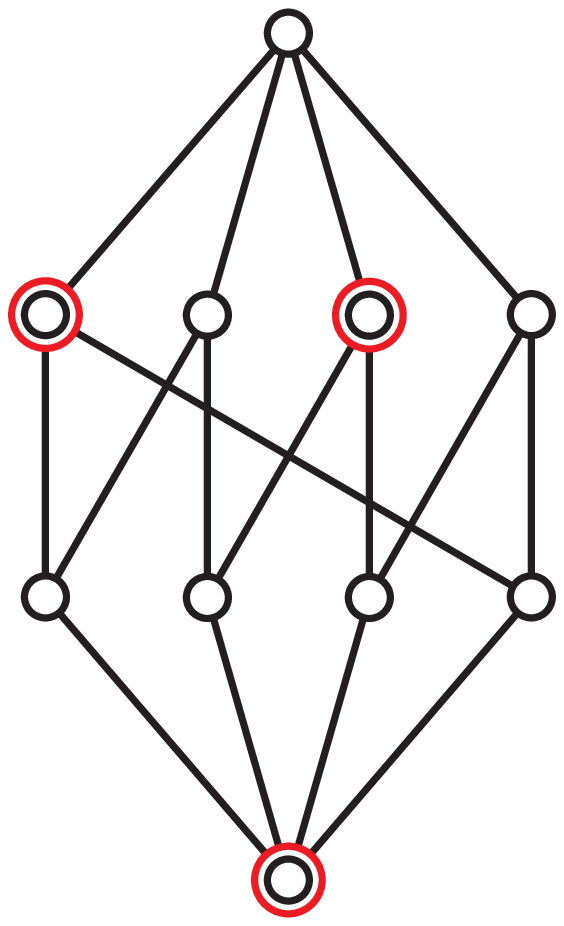 scaled 400}}
\rput(7.95,2.6){{\red $e_1$}}\rput(9.95,2.6){{\red $e_2$}}
\rput(9.95,.3){{\red $e_1\wedge e_2$}}
\psline[linewidth=.2mm]{<-}(10.4,1.45)(11,1.45)
\psline[linewidth=.2mm]{<-}(10.4,2.6)(11,2.6)
\rput(11.2,1.45){$U_i$}\rput(11.2,2.6){$X_i$}
}
\end{pspicture}
  \caption{The homomorphism $f$ of Theorem \ref{section:renner:result100} 
need not be injective.}
\label{section:renner:figure300} 
\end{figure}

Figure \ref{section:renner:figure300} 
gives the lattice of idempotents of the reflection monoid associated to
$\M$ (left) with a pair $\ve_1\ve_2\not=0$ marked, mapping via $f$
to $e_1\wedge e_2=0$ (right).
\end{example}

\begin{example}\label{section:renner:example300} 
Not only does the homomorphism $f$ fail to be injective in Example 
\ref{section:renner:example200},
but we can also show quite easily that $R_\M$ cannot be isomorphic to a
reflection monoid. For, suppose that $R_\M \cong M(W,\cS)$ where $\cS$
is a system of subspaces of a Euclidean space $V$ on which $W$ acts as a
reflection group. Since $W$ must be isomorphic to the group of units of
$R_\M$, we have $W = W(B_2)$. Hence four of the elements of order 2 in $W$
must be reflections. Also, the lattice $\cS$ must be isomorphic to the
lattice shown on the right in
Figure~\ref{section:renner:figure300}. Moreover, if the unique minimal 
element of $\cS$ is a non-zero subspace, we can factor it out to obtain
a lattice of subspaces with minimal element $\{0\}$.

Reading from left to right, let the elements of $\cS$ indicated 
in Figure \ref{section:renner:figure300} be
$U_0,U_1,U_2,U_3$ and $X_0,X_1,X_2,$ $X_3$ respectively. The intersection
of any two $U_i$'s is zero, as is the intersection of $X_0$ and $X_2$.
Hence for any choice of non-zero vectors $\uu_i \in U_i$ ($i =
0,1,2,3$), the set $\{ \uu_0^{},\dots,\uu_3^{} \}$ is linearly
independent.

The group of units of $R_\M$ is the automorphism group of $E(R_\M)$
where the action is by conjugation. Hence $W$ acting by conjugation  on
$\{ \ve_Y^{} \,|\, Y \in \cS\}$ gives all automorphisms of $E(M(W,\cS))$
and since $\ve_{Yg}^{} = g^{-1}\ve_Y^{}g$ for all $Y\in \cS$ and $g\in
W$, the $W$-action on $\cS$ gives all the automorphisms of $\cS$.

Now, automorphisms of $\cS$ are determined by their effect on the $U_i$. Let
$g,g'\in W$ be such that their actions give rise to the automorphisms
determined by interchanging $U_0$ with
$U_3$ and $U_1$ with $U_2$, and interchanging $U_0$ with
$U_1$ and $U_2$ with $U_3$ respectively.  Choose
$\uu_i\in U_i$ for $i=0,1$; then $\uu_0^{}g\in U_3$ and $\uu_1^{}g \in U_2$, so that 
$\{\uu_0^{},\uu_1^{},\uu_0^{}g,\uu_1^{}g\}$ is a basis for the subspace
it spans, say $U$. It is readily verified that $-1$ is an eigenvalue of
$g|_U^{}$ of multiplicity 2, so that $-1$ cannot be a simple eigenvalue
of $g$ itself. Thus $g$ (which has order 2) is not a reflection. 
Similarly, $g'\not= g$ is not a reflection. This is a contradiction since
there is only one element of order 2 in $W$ which is not a reflection. 
\end{example}





We conclude the subsection by mentioning that several authors have
calculated the orders of certain Renner monoids. The most general
results (which include all earlier ones) are in \cite{ZhenhengLi06}.
We will analyse in more detail the connection between reflection monoids and 
linear algebraic monoids in a future paper. 

\section{Complements}
\label{section:complements}

In this final section we elaborate on a number of miscellaneous
issues thrown up in earlier sections, but not strictly part
of the flow of the paper.

\subsection{Factorizable inverse monoids}
\label{section:complements:factorizable}

We first met factorizable inverse monoids in \S\ref{section:generalities} where we
characterized the factorizable inverse submonoids of $\ml(V)$ as being the 
monoids $M(G,\cS)$ that form the main characters in our story. 
The results of that section suggest that, in a suitably ``de-linearized'' form, 
they can be used to provide a description of \emph{all\/} factorizable inverse monoids.

We take our cue from group theory, where the ``Cayley'' representation embeds a group
$G$ in the symmetric group $\SS_G$. The equivalent for an inverse monoid $M$ is the
Vagner-Preston representation \cites{Howie95,Lawson98}, 
which is a faithful representation
$M\hookrightarrow\II_M$ given by partial right
multiplication. Thus any characterization of inverse monoids (up to isomorphism)
can be restricted to the inverse submonoids of the symmetric inverse monoid.

Throughout this section, let $X$ be an arbitrary set.
We observe that if $M$ is an inverse submonoid of $\II_X$, then
$E=E(M) = M \cap  E(\II_X) = \{ \ve_Y^{}\,|\, Y = \dom\, \aa \text{ for some } 
\aa \in M\}$.
Equally, $E = \{ \ve_Y^{}\,|\, Y = \im\, \aa \text{ for some } \aa \in
M\}$ since $\im \aa = \dom\, \aa^{-1}$ for all $\aa \in M$.
Putting
$$\cS = \{\, \dom\, \aa \,|\, \aa \in M\},$$
we see that $\cS$ is a meet semilattice isomorphic to $E$. Moreover,
$X\in \cS$ since $M$ is a submonoid, and finally, if $Y\in \cS$ and 
$g \in G=G(M)$, then $Yg = \im (\ve_{Y}^{}g) \in
\cS$. Thus $\cS$ provides an example of a \emph{system of subsets\/}
in $X$ for $G$: a collection $\cS\subset\mathbf{2}^X$ such that
$X\in\cS$, $\cS G=\cS$ and $X\cap Y\in\cS$ for all $X,Y\in\cS$. If
$G\subset\SS_X$ is a group and $\cS$ a system in $X$ for $G$ then we
form the \emph{monoid of partial permutations\/}
$M(G,\cS):=\{g_Y^{}\,|\,g\in G, Y\in\cS\}\subset\II_X$. 

Note that if $g_Y^{},h_Z^{} \in M(G,\cS)$, then $(g_Y^{})^{-1} =
(g^{-1})_{Yg}^{} \in M(G,\cS)$ and $g_Y^{}h_Z^{} = (gh)_T^{}$ 
with $T = Y \cap Zg^{-1}$, so
that $M(G,\cS)$ is an inverse submonoid of $\II_X$. Clearly, $G$ is the group
of units, and the idempotents are $E= \{ \ve_Y^{} \,|\, Y\in \cS\}$. 
Moreover, every
element is by definition a restriction of a unit, so $M(G,\cS)$ is factorizable.
Here is the promised characterization:

\begin{proposition}\label{section:complements:factorizable:result100}
$M$ is a factorizable inverse monoid if and only if there is a set $X$, a group
$G\subset\SS_X$, and a system $\cS$  in $X$ for $G$, with
$M$ isomorphic to $M(G,\cS)$.
\end{proposition}

We have already seen that monoids of the form $M(G,\cS)$
are factorizable inverse submonoids of $\II_X$. For the converse, it suffices, 
by the Vagner-Preston representation, to assume
$M\subset\II_X$ for some $X$. Let $G$ be its group of units, $\cS=\{\,
\dom\, \ss \,|\, \ss \in M\}$ 
the system above, and form $M(G,\cS)$. Now proceed as in the proof of
Proposition \ref{section:generalities:result100}.

We also have:

\begin{theorem}\label{section:complements:factorizable:result200}
Let $G\subset\SS_X$ be finite and $\cS$ a finite system in $X$ or
$G$. Then $|M(G,\cS)|=\sum_{Y\in\cS}[G:G_Y]$ with $G_Y=\{g\in
G\,|\,yg=y\text{ for all }y\in Y\}$. 
\end{theorem}

The proof is identical to Theorem \ref{section:reflectionmonoids:result100}.

\subsection{Fundamental inverse monoids}
\label{section:complements:fundamental}

We extend the themes of the previous section to describe
another abstract class of inverse monoids of interest: the
fundamental inverse monoids. On any inverse monoid $M$, define the
relation $\mu$ by the rule:
$$ a\,\mu\, b \text{ if and only if }  a^{-1}ea = b^{-1}eb \text{ for all
} e\in E.$$
It is easy to see that $\mu$ is a congruence on $M$; it is
idempotent-separating in the sense that distinct idempotents in $M$ are
not related by $\mu$, and, in fact, it is the greatest
idempotent-separating congruence on $M$. We say that $M$ is
\textit{fundamental} if $\mu$ is the equality relation; in general,
$M/\mu$ is fundamental. 

The Munn semigroup \cite{Howie95}*{\S5.4} $\TT_E$ of a semilattice $E$ is defined to be the 
set of all isomorphisms $Ee\rightarrow Ef$ where $e,f\in E$ with $Ee\cong Ef$.
We have $\TT_E$
an inverse submonoid of $\II_E$ whose semilattice of idempotents is
isomorphic to $E$ (see \cite{Howie95}*{Theorem 5.4.4} or \cite{Lawson98}*{Theorem
5.2.7}). 

Given any inverse monoid $M$ and $a\in M$, define an element $\aa_a
\in \TT_{E(M)}$ as follows. The domain of $\aa_a$ is $Eaa^{-1}$ and
$x\aa_a = a^{-1}xa$ for $x\in Eaa^{-1}$. Note that $\im \aa_a = Ea^{-1}a$.
The main results (see \cite{Howie95}*{Theorems
5.4.4 and 5.4.5} or \cite{Lawson98}*{Theorems 5.2.8 and 5.2.9}) are that
the mapping $\aa:M \to \TT_{E(M)}$
given by $a\aa = \aa_a$ is a homomorphism onto a full inverse
submonoid of $\TT_{E(M)}$ such that $a\aa = b\aa$ if and only if
$a\,\mu\,b$. Moreover, an inverse monoid $M$ is fundamental if and only if $M$ is isomorphic to
a full inverse submonoid of $\TT_{E(M)}$.

The homomorphism $\aa:M \to \TT_{E(M)}$ of is called the
\textit{fundamental} or \textit{Munn
representation} of $M$. Note that $M$ is fundamental if and only if
$\aa$ is one-one.
 
It is well known that the symmetric inverse monoid 
$\II_X$ is fundamental for any set $X$--see, for
example, \cite{Howie95}*{Chapter 5, Exercise~22}. In contrast, for any
nonempty set $X$, it is easy to see that the 
monoid of partial signed permutations $\JJ_X$ is
\textit{not} fundamental: a simple calculation shows that the identity
of $\JJ_X$ and the transposition $(x,-x)$ are $\mu$-related.  
In Proposition \ref{examples:permutationmonoidsresult200} we saw 
that $\JJ_n$ is a reflection monoid, so there certainly are
non-fundamental reflection monoids. 

We now describe fundamental factorizable inverse monoids  in terms of 
semilattices and their automorphism
groups. We
remark that the principal ideals of a semilattice $E$ regarded as a monoid are
precisely the principal order ideals of $E$ regarded as a partially
ordered set. It will be convenient to write $\ve_x$ for the partial
identity with domain $Ex$. 

\begin{proposition}\label{reflection_monoids:result50}
\begin{description}
\item[(i).] If $E$ is a semilattice with unique maximal element
and $G$ is a subgroup of the automorphism group $\aut(E)$,
then the collection 
$\cS = \{ Ex \,|\, x\in E\}$
of all principal
ideals of $E$ forms a system in $E$ for $G$,
and
the resulting $M(G,\cS)\cong\langle G,E\rangle\subset\TT_E$.
\item[(ii).] If $M$ is a fundamental factorizable inverse monoid with group of units $G$
and idempotents $E$ then $M\cong\langle G,E\rangle\subset\TT_E$.
\end{description}
\end{proposition} 

The principal example for us is the
Renner monoid of a reductive monoid with $0$ as in
\S\ref{section:renner}: it is fundamental factorizable
by \cite{Putcha88}*{Proposition 11.1} with units the Weyl group
$W_\G$ and idempotents $E=E(\ov{T})$, thus
$R_\M\cong\langle W_\G, E\rangle\subset\TT_E\subset\II_E$.

\begin{proof}
Given $E$ and $G$ we observe that 
$\cS$ does  form a system in $E$ for $G$ since $E=
E\hat{1}$, for $\hat{1}$ the maximal element,
 $Ex \cap Ey = Exy$ and the image under $g\in G$ of $Ex$ is
$E(xg)$. 
We can thus define 
the factorizable inverse monoid $M(G,\cS)\subset\II_E$ as above. As $G$
is a subgroup of $\aut(E)$, it is a subgroup of the group of units of
$\TT_E$, and hence if $\ve_xg \in M(G,\cS)$ with $g\in G$, then $\ve_xg
\in \TT_E$. Thus $M(G,\cS)\subset \TT_E$; in fact, it is clearly a full
inverse submonoid of $\TT_E$ and so it is fundamental. Identifying $E(\TT_E)$
with $E$, it is also
clear that $M(G,\cS)$ is generated as a submonoid by $G$ and $E$.

For part (ii), $M$ is isomorphic to a full submonoid of $\TT_E$ by the
injectivity of the Munn representation, and we identify this submonoid
with $M$. The group $G$ is a 
subgroup of the group
of units of $\TT_E$, that is, of $\aut(E)$. As above $\cS = \{\dom\,\aa
\,|\, \aa \in M\}$ is a system in $E$ for $G$, and since $M$
is factorizable we have $M = M(G,\cS)$. Thus $M$ is generated by $G$ and $E$
(identifying $E$ with $E(\TT_E)$).
\qed
\end{proof}

\parshape=15 0pt\hsize 0pt\hsize 0pt\hsize 0pt\hsize 0pt\hsize 
0pt\hsize 0pt\hsize 0pt\hsize 0pt.7\hsize 
0pt.7\hsize 0pt.7\hsize 0pt.7\hsize 0pt.7\hsize 0pt.7\hsize
0pt.7\hsize 
We finish by returning to reflection monoids
and giving an example of a non-fundamental
reflection monoid in which the restriction of the Munn
representation to the group of units is one-one. 
(We have seen that $\JJ_X$
is not fundamental, but in this case there are distinct units which 
are $\mu$-related.)
First,  note that if $M
= M(W,\cS)$ is a reflection monoid, and $\aa \in M$ has domain $X$, then
for any $Y \in \cS$ we have 
\begin{equation} \label{equation}
               \aa^{-1}\ve_Y^{}\aa = \ve_{(Y\,\cap\,X)\aa}^{}.
\end{equation} 
Now let $V=\mathbb{R}^2$ 
and $W$ the reflection group of either of the two triangles shown
(so $W \cong \SS_3$). 
The $\R$-spans of the vectors shown, together
with $V$ and $0$, form a system (of subspaces) for $W$. Let
$\rho\in W$ be the rotation through $2\pi/3$ and $\tau\in W$ the
reflection in the $y$-axis.
The $\mu$-class
of the identity $\ve_V^{}$ is a normal subgroup of $W$ and so to show that
$\mu$ is trivial on $W$, it is enough to show that $\rho$ and $\ve_V^{}$
are not $\mu$-related 
This is clear from \eqref{equation} using any of the six lines
for $Y$. 
\vadjust{\hfill\smash{\lower 0pt
\llap{
\begin{pspicture}(0,0)(3,3)
\rput(1.2,1.9){\BoxedEPSF{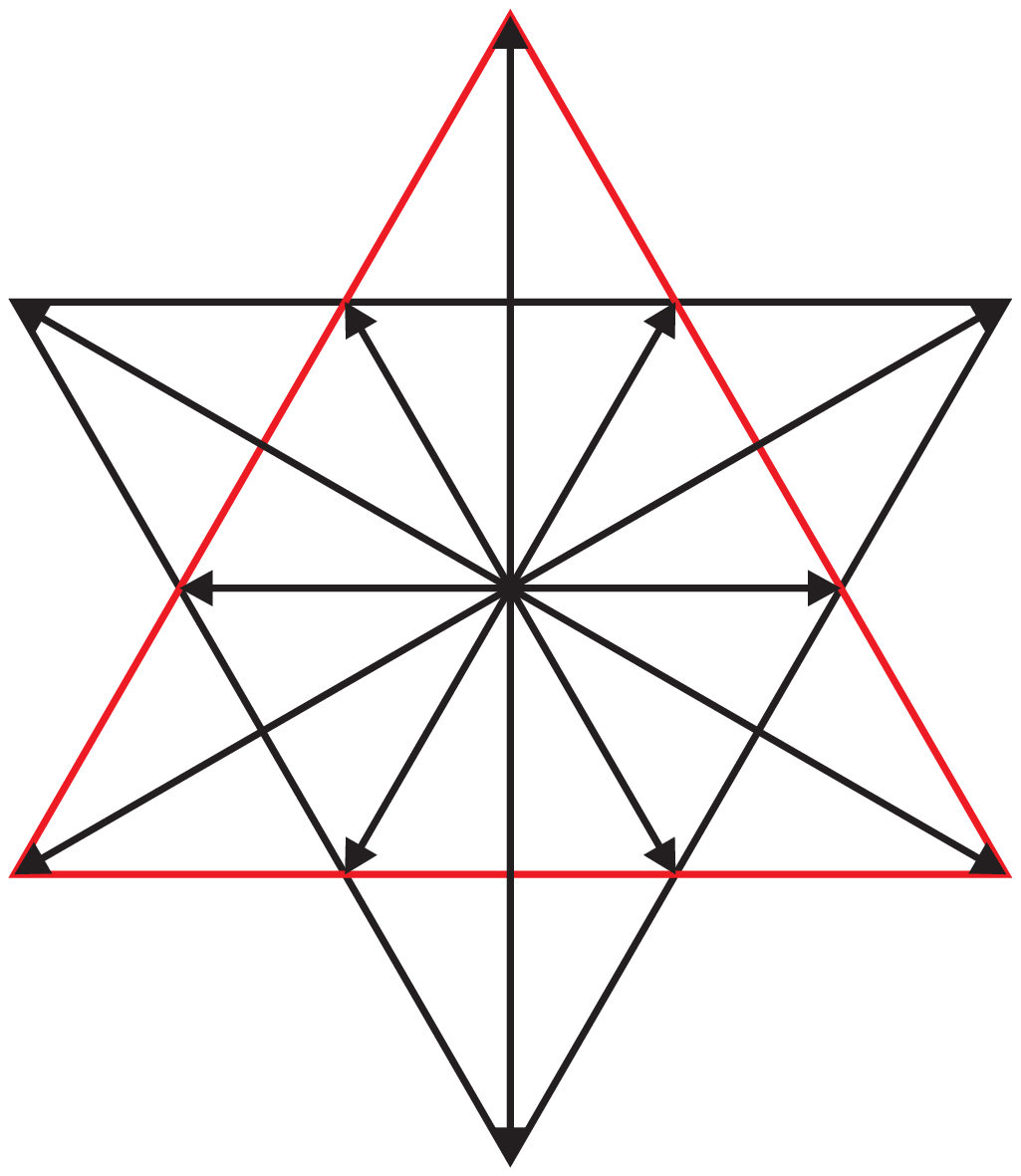 scaled 400}}
\end{pspicture}}}}\ignorespaces
On the other hand, letting 
$X$ be the $x$-axis, we see that $\tau_X^{}$ and $\ve_X^{}$ are
distinct but $\mu$-related.




\end{document}